\documentclass[preprint,10pt]{elsarticle}
\usepackage{amssymb}
\usepackage{amsfonts}
\usepackage{amsmath}
\usepackage{mathrsfs}
\usepackage{psfrag}
\usepackage{latexsym}
\usepackage{color}
\usepackage{amsbsy}
\usepackage{epsfig}
\usepackage{enumerate}
\usepackage{bm}

\usepackage{geometry}
\usepackage{graphicx}
\usepackage{float}
\usepackage{subfigure}
\usepackage[]{caption2} %新增调用的宏包
\renewcommand{\thesubfigure}{(\roman{subfigure})}% 此外，还可设置图编号显示格式，加括号或者不加括号
\makeatletter \renewcommand{\@thesubfigure}{\thesubfigure \space}%子图编号与名称的间隔设置
\renewcommand{\p@subfigure}{} \makeatother

\newtheorem{thm}{Theorem}[section]
\newtheorem{prop}{Proposition}[section]
\newtheorem{coro}{Corollary}[section]
\newtheorem{example}{Example}[section]
\newtheorem{lem}{Lemma}[section]
\newtheorem{definition}{Definition}[section]
\newtheorem{remark}{Remark}[section]

\newcommand{\Proof}{\textbf{Proof}}

\journal{ }

\begin{document}

\begin{frontmatter}
\title{On Independence for Capacities to Fit Ellsberg's Model with a Weak Law of Large Numbers}
\author[author1,author2]{Weihuan Huang}
\author[author3]{Yiwei Lin}\ead{lin$\_$yiwei@126.com}

\address[author1]{Zhongtai Securities Institute for Financial Studies, Shandong University, Jinan 250100, China}
\address[author2]{Department of Economics, The University of Kansas, Lawrence 66045, USA}
\address[author3]{School of Mathematics, Shandong University, Jinan 250100, China}

\cortext[cor]{The first author Huang's work has been supported by National Natural Science Foundation of China (Grant Nos. 11701331 and 11601280 and 11871050) and Natural Science Foundation of Shandong province (Grant No. ZR2017QA007), and Huang thanks the financial support from China Scholarship Council and Shandong University when staying at the University of Kansas.}

\begin{abstract}
This paper introduces new notions of Fubini independence and Exponential independence of random variables under capacities to fit Ellsberg's model, and finds out the relations between Fubini independence, Exponential independence, MacCheroni and Marinacci's independence and Peng's independence. As an application, we give a weak law of large numbers for capacities under Exponential independence. Simulations show that Ellsberg's model enjoy the weak law of large numbers when there is mean uncertainty with or without variance uncertainty.
\end{abstract}
\end{frontmatter}

\section{Introduction}

\iffalse
Ellsberg's model is a well-known model in economics with ambiguity. Consider a sequence of Ellsberg urns, you are told that each contains 100 balls, for example, that are either red or blue, thus $S=\{R,B\}$. You may also be given additional information, symmetric across urns, but it does not pin down either the precise composition of each urn or the relationship between urns.
\fi

Traditional philosophical wisdom on the notion INDEPENDENCE comes from the classical probability models with a finite sample space, e.g. draw red balls from different urns in which the propotion of red balls in each urns are known, and it is based on the fact that the probability (comes from the known propotion or the frequency of repeated experiment) is prior. With probability theory, we can measure risk by variance in traditional finance. However, Ellsberg \cite{Ellsberg} found that there are uncertainty that are not risk, i.e. in the above example we can neither know the propotion of red balls prior nor do experiments to find the exact probability by frequency. This is the famous Ellsberg's Paradox (model) in modern finance and economics which leads us to find the INDEPENDENCE and behavior of the frequency (law of large numbers) under uncertainty. In this paper, we will find a Fubini independence and Exponential independence to fit Ellsberg's model, and a weak law of large numbers under these INDEPENDENCE.

Peng \cite{Peng2006} has developed the non-linear expectations theory to model uncertainty (ambiguity) and coherent risk measures. In this theory, Peng has proved a new law of large numbers and a central limit theorem under sublinear expectations with a new INDEPENDENCE condition, which laid the theoretical foundations for the non-linear expectation framework. The strength of non-linear expectation theory is utilizing a kind of nonlinear heat equations to construct $G$-Brownian Motion ($G$-normal distribution) and maximal distribution with uncertainty. However, a weakness of the theory is the non-linear expectation can not measure indicative functions, i.e. can not use a non-additive probability (capacity) to measure random events. This leads to failure to fit the Ellsberg's model. Following example shows that we can not use Peng's independence for upper-expectation to model Ellsberg's urns.

\begin{example}
	Peng's independence: a r.v. $Y$ is said to be Peng's independent from another r.v. $X$ under sub-linear expectation $$\mathbb{E}[\cdot]:=\sup_{p\in\mathcal{P}}E_P[\cdot],$$ where $\mathcal{P}$ is a set of probabilities, if for each test function $\phi\in C_{l.Lip}(\mathbb{R}^2)$ we have
	\begin{equation}\label{Pengindep}
	\mathbb{E}[\phi(X,Y)] = \mathbb{E}[\mathbb{E}[\phi(x,Y)]|_{x=X}].
	\end{equation}
	In this paper, we call \eqref{Pengindep} Peng's independence.
	
	We consider a special Ellsberg's urns (two urns) for example. Let $\Omega_i=\{R,B\}$, $i=1,2$ be two urns (i.e. each urn has red balls and black balls). And two r.v.s $X$ and $Y$, for $\omega_1\in\Omega_1$, $\omega_2\in\Omega_2$
	\begin{eqnarray*}X(\omega_1)=
		\begin{cases}
			1, &\omega_1=R ,\cr 0, &\omega_1=B,
		\end{cases}
	\end{eqnarray*}
	\begin{eqnarray*}Y(\omega_2)=
	\begin{cases}
		1, &\omega_2=R ,\cr 0, &\omega_2=B,
	\end{cases}
\end{eqnarray*}
and a function $\phi\in C_{l.Lip}(\mathbb{R}^2)$
	\begin{eqnarray*}\phi(0,y)=
	\begin{cases}
		1, &y=1 ,\cr 0, &y=0,
	\end{cases}
\end{eqnarray*}
	\begin{eqnarray*}\phi(1,y)=
	\begin{cases}
		0, &y=1 ,\cr 1, &y=0,
	\end{cases}
\end{eqnarray*}
 let $\sup\limits_{P\in\mathcal{P}}P(\omega_1=R)=\sup\limits_{P\in\mathcal{P}}P(\omega_2=R)=0.5$, $\inf\limits_{P\in\mathcal{P}}P(\omega_1=R)=\inf\limits_{P\in\mathcal{P}}P(\omega_2=R)=0.3$ and $\mathcal{P}:=\{ P\ {\rm is\ a\ prob.} | P\leq \sup\limits_{P\in\mathcal{P}}P\}$, which is a core. In this case, it is easy to check that $\mathbb{E}[\phi(X,Y)] \not= \mathbb{E}[\mathbb{E}[\phi(x,Y)]|_{x=X}]$, in fact $\mathbb{E}[\phi(X,Y)]=0.5$ and $\mathbb{E}[\mathbb{E}[\phi(x,Y)]|_{x=X}]=0.6$, which means that Peng's independent do not fit the Ellsberg's model for all local Lipschitz functions even though $\omega_1$ and $\omega_2$ are from two different urns.
\end{example}

Ghirardato \cite{Ghi} has proved the Fubini theorem for capacities, and the key property is that the test functions are slice-comonotonicity (see Definition \ref{slicecomonotonic}). In this paper, we restrict the test functions $\phi$ in Peng's independence \eqref{Pengindep} to be slice-comonotonic to fit the Ellsberg's model, see section \ref{example}, and we called it Fubini independence, see Definition \ref{Fubiniindep}. And we find out that Fubini independence can implies Exponential independence, see Definition \ref{Expindep}, which also fits Ellsberg's model. Meanwhile, we prove a weak law of large numbers for capacities under Exponential independence to describe the behavior of frequency.

In capacity theory, MacCheroni and Marinacci \cite{MM} gives a strong law of large numbers for capacities under the independence condition as follows. $\{X_n\}_{n\geq 1}$ of r.v.s are pairwise independent with respect to capacities $\nu$ if, for each $n,m\geq 1$ and for all open subsets $G_n$, $G_m$ of $\mathbb{R}$,
$$
\nu(\{X_n\in G_n, X_m\in G_m\}) = \nu(\{X_n\in G_n\})\nu(\{X_m\in G_m\}).
$$
We call it MacCheroni and Marinacci's independence in our paper. Obviously, this independence fits Ellsberg's model too, but under this condition they should use more condition on $X_n$ or $\nu$ to prove law of large numbers, i.e. $X_n$ are bounded continuous or $\nu$ is continuous. We do not use these strong conditions to prove law of large numbers under Exponential independence.

The paper proceeds as follows. In section 2, we give definitions of Fubini independence and Exponential independence under capacities. In section 3, we give the relations between Fubini independence, Exponential independence and MacCheroni and Marinacci's independence. In order to find out that Fubini independence implies Exponential independence, we prove a Fubini Theorem. In section 4, we explain why Fubini independence fit the Ellsberg's model. In section 5, we prove a weak law of large numbers for capacities under Exponential independence as an application. In section 6, we show that Ellsberg's model enjoys the weak law of large numbers by simulations. We put the proof of Theorem \ref{fubiniduli} in appendix.

\section{Preliminaries}
\subsection{Basic concepts and lemmas}

Let $\Omega$ be a sample space and $\mathcal{F}$ be its $\sigma$-algebra. A set function $V:\Omega\rightarrow [0,1]$ is called a capacity if it satisfies the following:

(i) $V(\emptyset)=0,\ V(\Omega)=1$;

(ii) $\forall A,\ B\in \mathcal{F},\ A\subseteq B$, $V(A)\leq V(B)$.

A capacity is called total monotonicity if it satisfies the additional property:

(iii)(total monotonicity) $\forall n>0,$ and every collection $A_1,\cdots,A_n\in \mathcal{F}$, $$V(\cup_{i=1}^{n}A_i)\geq \sum_{\emptyset\neq I\subseteq \{1,\cdots,n\}}(-1)^{|I|+1}V(\cap_{i\in I}A_i).$$

A capacity is called total alternating if it satisfies the additional property:

(iv)(total alternating) $\forall n>0,$ and every collection $A_1,\cdots,A_n\in \mathcal{F}$, $$V(\cup_{i=1}^{n}A_i)\leq \sum_{\emptyset\neq I\subseteq \{1,\cdots,n\}}(-1)^{|I|+1}V(\cap_{i\in I}A_i).$$

We call it 2-alternating, if (iv) holds only for $n=2$.

Let $(\mathbb{R},\mathcal{B})$ denotes real numbers space $\mathbb{R}$ with the set $\mathcal{B}$ of all its Borel sets, $X_1, X_2, \cdots, X_n$ are random variables (r.v.) on $(\Omega, \mathcal{F}, V)$, i.e. $X_i:\Omega\rightarrow \mathbb{R}$, $i=1, 2, \cdots, n,$ are $\mathcal{F}$-measurable. The Choquet integral of a bounded r.v. X with respect to any capacity $V$ is defined by
$$
E_V[X]:=\int_{0}^{+\infty}V(\{X\geq t\})dt+\int_{-\infty}^0[V({X\geq t})-1]dt,
$$
where the integrals on the right hand are Riemann integrals. In this paper, $E_V$ denotes the Choquet integral with capacity $V$, $E_P$ denotes classical expectation with probability $P$, $\mathbb{E}_V$ denotes the upper envelop with respect to the core of capacity $V:=\sup\limits_{P\in\mathcal{P}}P$, i.e.
$$
\mathbb{E}_V[X]:= \sup_{P\in\mathcal{P}}E_P[X],
$$
where core means that $\mathcal{P}=\{P\ {\rm is\ a\ prob.}| P\leq V\}$. If $V$ is a 2-alternating capacity, we denote it by $\mathbb{V}$.

We define a product space $(\mathbb{R}^n,\mathcal{B}^n)$, where $\mathcal{B}^n$ is the product Borel sets on $\mathbb{R}^n$, that is, the smallest $\sigma$-algebra of subsets of $\mathbb{R}^n$ which contains all rectangles, and let $\mathcal{R}$ be the set of all rectangles on $\mathbb{R}^n$.
Any capacity $\sigma$ on $(\mathbb{R}^2,\mathcal{B}^2)$ will be defined a product capacity, its marginals on $\mathbb{R}$ will be respectively the capacities $\mu$ and $\nu$ as follows, for all $A, B\in\mathcal{B}$,
$$
\mu(A):=\sigma(A\times\mathbb{R}), \ \ \nu(B):=\sigma(\mathbb{R}\times B).
$$

\iffalse
\begin{definition}[Fubini Independent Product Capacity]
A capacity $\sigma$ on $(\mathbb{R}^2, \mathcal{B}^2)$ with marginals $\mu$ and $\nu$ is said to be Fubini Independent for $\mathcal{H}$ (for $\mathcal{B}^2$ or for $\mathcal{R}$), if for every set $A\in\mathcal{H}$ ($\mathcal{B}^2$ or $\mathcal{R}$), we have
\begin{equation}
\sigma(A)=\mathbb{E}_{\nu} [\mu((X,y)\in A)|_{y=Y}].
\end{equation}
\end{definition}
\fi

%$X_i$ deduces a capacity $V_{X_{i}^{-1}}$ on $\mathbb{R}$ as $V_{X_{i}^{-1}}(A):=V(X_i\in A)$, $\forall A\in \mathcal{B}$.

Let $\mathcal{U}_{f}$ be the set of all upper intervals in $\mathbb{R}^n$ of the form $\{(x_1, \cdots, x_n)\in\mathbb{R}^n : f(x_1, \cdots, x_n)>\alpha\}$, for some $\alpha\in\mathbb{R}$, i.e. %or $\{(x_1,x_2)\in\mathbb{R}^2 : f(x_1,x_2)>\alpha\}$,
$$
\mathcal{U}_{f}=\{  \{(x_1, \cdots, x_n)\in\mathbb{R}^n : f(x_1, \cdots, x_n)>\alpha\} | \alpha\in\mathbb{R} \}.
$$
Obviously, $\mathcal{U}_{f}$ is a chain (a family completely ordered by set inclusion).

\begin{definition}\label{comonotonic}
Let $f_1, f_2:\mathbb{R}^2\rightarrow\mathbb{R}$, we say that $f_1$ and $f_2$ are comonotonic, if for every $(x_1,x_2), (x_1',x_2')\in\mathbb{R}^2$, we have
$$
[f_1(x_1,x_2)-f_1(x_1',x_2')][f_2(x_1,x_2)-f_2(x_1',x_2')]\geq 0,
$$
a class of functions $\mathcal{G}$ is said to be comonotonic if for every $f_1, f_2\in\mathcal{G}$, $f_1$ and $f_2$ are comonotonic.
\end{definition}

\begin{definition}\label{slicecomonotonic}
$f:\mathbb{R}^2\rightarrow\mathbb{R}$, $f$ has comonotonic $x_1$-sections, if for every $x,x'\in\mathbb{R}$, $f(x,\cdot)$ and $f(x',\cdot)$ are comonotonic functions. If $f$ both has comonotonic $x_1$- and $x_2$- sections, then we call it is slice-comonotonic. A set $A\in\mathbb{R}^2$ is said to be comonotonic, if its characteristic function has comonotonic $x_1$-sections.
\end{definition}
Similarly, we can define comonotonic $x_i$-sections in $\mathbb{R}^n$, and comonotonic functions in $\mathbb{R}^n$.

For
 \begin{equation}\label{fhat}
 \hat{f}(x_1,x_2):=e^{\varphi_1(x_1)+\varphi_2(x_2)},
 \end{equation}
we have
\begin{align*}
&[\hat{f}(x_1,x_2)-\hat{f}(x_1,x_2')][\hat{f}(x_1',x_2)-\hat{f}(x_1',x_2')]\\
=& [e^{\varphi_1(x_1)+\varphi_2(x_2)}-e^{\varphi_1(x_1)+\varphi_2(x_2')}]\cdot[e^{\varphi_1(x_1')+\varphi_2(x_2)}-e^{\varphi_1(x_1')+\varphi_2(x_2')}] \\
=&e^{\varphi_1(x_1)}e^{\varphi_1(x_1')}[e^{\varphi_2(x_2)}-e^{\varphi_2(x_2')}]^2 \geq 0,
\end{align*}
hence, $\hat{f}(x_1,x_2)$ has comonotonic $x_1$-sections, similarly, it has comonotonic $x_2$-sections, i.e. $\hat{f}(x_1,x_2)$ is slice-comonotonic. Similarly, $\hat{f}(x_1, \cdots, x_n):=e^{\varphi_1(x_1)+\cdots+\varphi_2(x_n)}$ is slice-comonotonic.

The following two lemmas will be used to prove the main results of our paper.

\begin{lem}[Lemma 1 in \cite{Ghi}]\label{Lem1}
Suppose that $\mathcal{G}$ is a comonotonic class of bounded r.v. from $\Omega$ into $\mathbb{R}$ and $V$ is a capacity on $(\Omega, \mathcal{F})$. Then we can find a probability measure $P$ on $(\Omega, \mathcal{F})$ such that for every $X\in\mathcal{G}$,
$$
E_{V}[X] = E_{P}[X].
$$
\end{lem}

\begin{lem}[Lemma 2 in \cite{Ghi}]\label{Lem2}
Let $f: \mathbb{R}^2 \rightarrow \mathbb{R}$ be a bounded, measurable function with comonotonic $x_1$- or $x_2$- sections. Then every $A\in\mathcal{U}_f$ is a comonotonic set.
\end{lem}
\begin{remark} $\hat{f}$ given by \eqref{fhat} are slice-comonotonic, so every $A\in\mathcal{U}_{\hat{f}}$ is a comonotonic set.
%\section{Fubini Theorem for Capacities}
\end{remark}

\subsection{Fubini Independence and Exponential Independence}

We introduce two new independence in this section to fit Ellsberg's model, the Fubini independence is inspired by Peng's independence in Peng \cite{Peng} and Fubini Theorem in Ghirardato \cite{Ghi}, and Exponential independence was also introduced to prove a strong law of large number in Chen, Huang and Wu \cite{CHW}.
\begin{definition}[Fubini Independence for Random Variables]
$X_n$ is said to be Fubini independent of $X_1, \cdots, X_{n-1}$ for $\mathcal{U}_{f}$ (for $\mathcal{B}^2$ or for $\mathcal{R}$),
if for every set $A\in\mathcal{U}_{f}$ ($\in\mathcal{B}^2$ or $\in\mathcal{R}$), we have
\begin{equation}\label{Fubini Independent}
V((X_1, \cdots, X_n)\in A)=E_{V}\left[V\left((x_1, \cdots, x_{n-1}, X_n)\in A\right)|_{x_1=X_1, \cdots, x_{n-1}=X_{n-1}}\right].
\end{equation}\label{Fubiniindep}
\end{definition}
Obviously, $X_1$ and $X_2$ are Fubini independent for $\mathcal{R}$ is equivalent to
\begin{equation}\label{MMindep}
V(X_1\in A, X_2\in B)=V(X_1\in A)V(X_2\in B),
\end{equation}
for all $A, B\in\mathcal{B}$, which is introduced by MacCheroni and Marinacci \cite{MM} to prove a strong law of large numbers for capacities, we call \eqref{MMindep} MacCheroni and Marinacci's independence.

If we consider $X_n$ is Fubini independent of $X_1, \cdots, X_{n-1}$ for $\mathcal{U}_{\hat{f}}$, where $\hat{f}$ is defined by \eqref{fhat}, we have (\ref{Fubini Independent}) as
\begin{align*}
V(&\varphi_1(X_1)+\cdots+\varphi_n(X_n)> \alpha) \\
&=E_{V}\left[V\left( \varphi_1(x_1)+\cdots+\varphi_{n-1}(x_{n-1})+\varphi_n(X_n)> \alpha\right)|_{x_1=X_1, \cdots, x_{n-1}=X_{n-1}}\right],
\end{align*}
for all $\alpha\in\mathbb{R}$.

\begin{definition}[Exponential Independence]\label{Expindep}
We called $X_n$ is exponential independent of $X_1, \cdots, X_{n-1}$, if
\begin{equation}\label{EI}
E_{V}\left[e^{\sum_{i=1}^{n}\varphi_{i}(X_i)}\right]=E_{V}\left[e^{\sum_{i=1}^{n-1}\varphi_{i}(X_i)}\right] \cdot E_{V}\left[e^{\varphi_{n}(X_n)}\right],
\end{equation}
where $\varphi_{i}, i=1, \cdots, n$ are bounded functions.
\end{definition}

\section{Relations between Independences}

In classical probability space $(\Omega,\mathcal{F},P)$, where $P$ is a probability measure, and it is easy to have following properties.

\begin{prop}\label{classical}
The following three are equivalent:

(I) For all Borel set $A\in\mathcal{B}^2$, we have $$P((X_1,X_2)\in A)=E_{P}[P\left((x_1,X_2)\in A\right)|_{x_1=X_1}];$$

(II) For all Borel set $A,B\in\mathcal{B}$, we have $$P(X_1\in A, X_2\in B)=P(X_1\in A)P(X_2\in B);$$

(III) For any bounded continuous functions $\varphi_1$ and $\varphi_2$, we have
$$
E_{P}\left[e^{\varphi_{1}(X_1)+\varphi_{2}(X_2)}\right]=E_{P}\left[e^{\varphi_{1}(X_1)}\right] \cdot E_{P}\left[e^{\varphi_{2}(X_2)}\right].
$$
\end{prop}

However, the cases is different under $(\Omega,\mathcal{F},V)$, where $V$ is a capacity.

\begin{prop}
(I) If random variables $X_1$ and $X_2$ are Fubini independent for $\mathcal{B}^2$, then $X_1$ and $X_2$ are Fubini independent for $\mathcal{U}_{f}$ and for $\mathcal{R}$; \\
(II) If random variables $X_1$ and $X_2$ are Fubini independent for $\mathcal{B}^2$, then for all $A, B\in\mathcal{B}$, we have
$$
V(X_1\in A, X_2\in B)=V(X_1\in A)V(X_2\in B);
$$
%(III) If random variables $X_1$ and $X_2$ are Fubini independent for $\mathcal{B}^2$, then for $a, b\in\mathbb{R}$, we have
%$$
%V(X_1\leq a, X_2\leq b)=V(X_1\leq a)V(X_2\leq b);
%$$
\end{prop}

\Proof. (I) is obvious.
(II), $V((X_1,X_2)\in A\times B)=V(X_1\in A,X_2\in B)$, and $E_V[V((X_1,x)\in A\times B)|_{x=X_2}]=E_V[E_V[I_A(X_1)]I_B(X_2)]=V(X_1\in A)V(X_2\in B)$.
%And (III) follows (II).

\begin{prop}
If $X_1$ and $X_2$ are exponential independent random variables, then for $a, b\in\mathbb{R}$, we have
$$
V(X_1\leq a, X_2\leq b)=V(X_1\leq a)V(X_2\leq b);
$$
and, for any Borel sets $A$ and $B$,
$$
V(X_1\in A, X_2\in B)=V(X_1\in A)V(X_2\in B).
$$
\end{prop}

\Proof. Let
\begin{eqnarray*}\varphi_{n,1}(x)=
\begin{cases}
0, &x\leq a,\cr -n, &x> a,
\end{cases}
\end{eqnarray*}
and
\begin{eqnarray*}\varphi_{n,2}(x)=
\begin{cases}
0, &x\leq b,\cr -n, &x> b,
\end{cases}
\end{eqnarray*}
are bounded functions, then by Dominated (Bounded) Convergence Theorem (Theorem 8.9 in Denneberg \cite{Den}),
\begin{align*}
E_{V}\left[e^{\varphi_{n,1}(X_1)+\varphi_{n,2}(X_2)}\right] &=E_{V}\left[\left(I_{\{X_1\leq a\}}+e^{-n}I_{\{X_1>a\}}\right)\cdot\left(I_{\{X_{2}\leq b\}}+e^{-n}I_{\{X_2>b\}}\right)\right] \\
&=E_V\left[I_{\{X_1\leq a, X_2\leq b\}}
+e^{-n}[I_{\{X_1\leq a, X_2>b\}}+I_{\{X_1>a, X_2\leq b\}}]
+e^{-2n}I_{\{X_1>a,X_2>b\}}\right] \\
&\rightarrow_{n\rightarrow\infty} E_V\left[I_{\{X_1\leq a, X_2\leq b\}}\right]=V(X_1\leq a, X_2\leq b),
\end{align*}
On the other hand,
\begin{align*}
E_{V}\left[e^{\varphi_{n,1}(X_1)+\varphi_{n,2}(X_2)}\right] &=E_{V}\left[e^{\varphi_{n,1}(X_1)}\right]\mathbb{E}_{V}\left[e^{\varphi_{n,2}(X_2)}\right] \\
&=E_V[I_{\{X_1\leq a\}}+e^{-n}I_{\{X_1>a\}}]\cdot E_V[I_{\{X_2\leq b\}}+e^{-n}I_{\{X_2>b\}}] \\
&\rightarrow_{n\rightarrow\infty} E_V[I_{\{X_1\leq a\}}]\cdot E_V[I_{\{X_2\leq b\}}]= V(X_1\leq a)V(X_2\leq b),
\end{align*}
the proof of the second result is similar as we can let
\begin{eqnarray*}\varphi_{n,1}(x)=
\begin{cases}
0, &x\in A,\cr -n, &x\notin A,
\end{cases}
\end{eqnarray*}
and
\begin{eqnarray*}\varphi_{n,2}(x)=
\begin{cases}
0, &x\in B,\cr -n, &x\notin B,
\end{cases}
\end{eqnarray*}
so we completed the proof.

\begin{coro}
If $X_n$ is exponential independent of $X_1, \cdots, X_{n-1}$, for all $n\geq 1$, i.e.
$$
E_{V}\left[e^{\sum_{i=1}^{n}\varphi_{i}(X_i)}\right]=\prod_{i=1}^n E_{V}\left[e^{\varphi_{i}(X_i)}\right],
$$
then for all Borel sets $A_i$, we have
$$
V\Big(\bigcap_{i=1}^n X_i\in A_i\Big)=\prod_{i=1}^n V(X_i\in A_i).
$$
\end{coro}

\begin{thm}\label{Fubini}
If $X_n$ is Fubini independent of $X_1, \cdots, X_{n-1}$ for $\mathcal{U}_{\hat{f}}$, where $\hat{f}(x_1,\cdots,x_n):=e^{\sum_{i=1}^{n}\varphi_i(x_i)}$, $\varphi_i$ are bounded, then $X_n$ is exponential independent of $X_1, \cdots, X_{n-1}$.
\end{thm}

Theorem \ref{Fubini} follows Fubini Theorem for capacities. We give Fubini Theorem for $n=2$, for $n>2$ cases is same.

\begin{thm}\label{fubiniduli}
 For functions $\hat{f}(x_1,x_2)=e^{\varphi_{1}(x_1)+\varphi_{2}(x_2)},\ \forall x_1, x_2 \in \mathbb{R}$, defined on the product space $\mathbb{R}^2$, where $\varphi_i, i=1,2$ are bounded functions, $X_1$ and $X_2$ are Fubini independent random variables for $\mathcal{U}_{\hat{f}}$, then we have
 \begin{equation}\label{Fubinijieguo}
 E_{V}[\hat{f}(X_1,X_2)]=E_{V}[E_{V}[\hat{f}(x_1,X_2)]|_{x_1=X_1}].
 \end{equation}
\end{thm}

The proof of Theorem \ref{fubiniduli} is similar to the proof of Fubini Theorem in \cite{Ghi}, so we put it in Appendix. By Theorem \ref{fubiniduli}, we can prove Theorem \ref{Fubini}.

\iffalse
From the proof, we can see that
\begin{coro} For functions $\hat{f}(x_1,x_2)=e^{\varphi_{1}(x_1)+\varphi_{2}(x_2)},\ \forall x_1, x_2 \in \mathbb{R}$, defined on the product space $\mathbb{R}^2$, where $\varphi_i, i=1,2$ are bounded functions, then the following two are equivalent:

(I) $X_1$ and $X_2$ are Fubini independent random variables for $\mathcal{U}_{f}$;

(II) $\forall g\in\{f\}\cup \{I_{A}|A\in\mathcal{U}_{f}\}$, we have
\begin{equation}
 \mathbb{E}_{V}[g(X_1,X_2)]=\mathbb{E}_{V}[\mathbb{E}_{V}[g(X_1,x)]|_{x=X_2}].
 \end{equation}
\end{coro}
\fi

\Proof. [The proof of Theorem \ref{Fubini}] Because of Fubini Theorem, we have
\iffalse
\begin{align*}
\mathbb{E}_{V}\left[e^{\varphi_{1}(X_1)+\varphi_{2}(X_2)}\right] &=\mathbb{E}_{V_{(X_1,X_2)^{-1}}}\left[e^{\varphi_{1}(x_1)+\varphi_{2}(x_2)}\right] \\
&=\mathbb{E}_{V_{X_2^{-1}}}\left[\mathbb{E}_{V_{X_1^{-1}}}[e^{\varphi_{1}(x_1)+\varphi_{2}(x_2)}]\right] \\
&=\mathbb{E}_{V_{X_1^{-1}}}\left[e^{\varphi_{1}(x_1)}\right]\cdot\mathbb{E}_{V_{X_2^{-1}}}\left[e^{\varphi_{2}(x_2)}\right] \\
&=\mathbb{E}_{V}\left[e^{\varphi_{1}(X_1)}\right]\cdot\mathbb{E}_{V}\left[e^{\varphi_{2}(X_2)}\right].
\end{align*}
\fi
\begin{align*}
E_{V}\left[e^{\varphi_{1}(X_1)+\varphi_{2}(X_2)}\right] &=E_{V}\left[E_{V}\left[e^{\varphi_{1}(x_1)+\varphi_{2}(X_2)}\right]\Big|_{x_1=X_1}\right] \\
&=E_{V}\left[E_{V}\left[e^{\varphi_{2}(X_2)}\right]e^{\varphi_{1}(X_1)}\right] \\
&=E_{V}\left[e^{\varphi_{1}(X_1)}\right] E_{V}\left[e^{\varphi_{2}(X_2)}\right].
\end{align*}

\iffalse
\begin{remark}
Theorem \ref{Fubini} shows that Fubini independence for $\mathcal{U}_f$ is a sufficient condition of exponential independence, meanwhile, Fubini independence for $\mathcal{U}_f$ is a necessary condition of exponential independence in the following sense: Suppose not, then there is a set $A\in\mathcal{U}_f$ s.t. (\ref{Fubiniduli}) does not hold, i.e.
$$V((X_1,X_2)\in A)\neq\mathbb{E}_V[V((X_1,x)\in A)|_{x=X_2}]=\lambda(A),$$
say $<$, i.e.
$$V((X_1,X_2)\in A)<\mathbb{E}_V[V((X_1,x)\in A)|_{x=X_2}]=\lambda(A),$$
hence we have
$$\mathbb{E}_V[f_p(X_1,X_2)]<\mathbb{E}_V[\mathbb{E}_V[f_p(X_1,x)]|_{x=X_2}],$$
for infinitely many $p\geq 1$, it contradicts the (\ref{shizi3}), so
$$\mathbb{E}_V[f(X_1,X_2)]<\mathbb{E}_V[\mathbb{E}_V[f(X_1,x)]|_{x=X_2}],$$
and then it contradicts the (\ref{Fubinijieguo}).
%It shows us that we can not use
%$$V(X_1\in A, X_2\in B)=V(X_1\in A)V(X_2\in B)$$
%to prove Fubini Theorem and Exponential Independence.
\end{remark}
\fi

\section{Fubini Independence fits Ellsberg's model}\label{example}

If $V$ is a 2-alternating capacity, we denote it by $\mathbb{V}$. By Proposition 3 in \cite{Sch}, we can restate the Fubini independent for $\mathcal{U}_{\hat{f}}$ by the upper envelop $\mathbb{E}$ as:
$$
\mathbb{V}((X_1, \cdots, X_n)\in A)=\mathbb{E}_{\mathbb{V}}\left[\mathbb{V}\left((x_1, \cdots, x_{n-1}, X_n)\in A\right)|_{x_1=X_1, \cdots, x_{n-1}=X_{n-1}}\right],
$$
for $A\in\mathcal{U}_{\hat{f}}$, i.e.
\begin{align*}
\mathbb{V}(&\varphi_1(X_1)+\cdots+\varphi_n(X_n)> \alpha) \\
&=\mathbb{E}\left[\mathbb{V}\left( \varphi_1(x_1)+\cdots+\varphi_{n-1}(x_{n-1})+\varphi_n(X_n)> \alpha\right)|_{x_1=X_1, \cdots, x_{n-1}=X_{n-1}}\right],
\end{align*}
for all $\alpha\in\mathbb{R}$.

\subsection{Finite sample case}
As an explanatory example for Fubini independence, consider a sequence of Ellsberg's urns. Next example shows that the Ellsberg's model (with finite sample space) satisfies the Fubini independence for $\mathcal{U}_{\hat{f}}$ under upper envelop $\mathbb{E}$.
%and in the example we consider the upper envelop $\mathbb{E}_{\mathbb{V}}$, i.e. $\mathbb{V}$ is a 2-alternating capacity.

\begin{example}\label{ellsberg'smodel}
Let $i\geq 1$, $\Omega_i:=\{\omega_1^i, \omega_2^i, \cdots, \omega_n^i\}$, $\mathcal{F}_i:=2^{\Omega_i}$, where $\Omega_i$ is stand for the $i$-th urn, $X_i(\omega_k^i):=x_k$, for $k=1,2,\cdots,n$, and $x_1\leq x_2\leq\cdots\leq x_n$, if not we just reorder elements of $\Omega_i$.
Denote $P^l(X_i=x_j):=p_{lj}$, for $l=1,2,\cdots,m$, $j=1,2,\cdots,n$, where $p_{lj}\geq 0$, $\sum_{j=1}^{n}p_{lj}=1$, and $\mathcal{P}_i:=\{P^l | l=1,2,\cdots,m\}$. \footnote{$p_{lj}\geq 0$ are arbitrary fixed, it means that $\mathcal{P}_i$ can be arbitrary. So the example is considered as general case for finite sample Ellsberg's model.}

We define a probability $P'_i$ on $\mathcal{F}_i$ as: $P'_i(X_i=x_n):=\max\limits_{1\leq l\leq m}p_{ln}$, and
$$
P'_i(X_i=x_k):=\max_{1\leq l\leq m} \sum_{j=k}^{n}p_{lj}-\max_{1\leq l\leq m} \sum_{j=k+1}^{n}p_{lj},
$$
where $k=1,2,\cdots,n-1$.

Obviously, $P'_i(X_i=x_k)$ is between $0$ and $1$ and we let it to be additive and monotone on $\mathcal{F}_i$, i.e.
$$
P'_i(X_i\in A):=\sum_{k=1}^{n} P'_i(X_i=x_k) \delta_{x_k}(A),
$$
for $A\in\mathcal{F}_i$, hence $P'_i$ is a probability on $\mathcal{F}_i$.

Consider the upper capacity $\mathbb{V}$ on $\mathcal{F}_i$ as
$$\mathbb{V}(X_i\in A):=\max_{P\in\mathcal{P}_i}P(X_i\in A),$$
for $A\in\mathcal{F}_i$, and upper envelop
$$
\mathbb{E}[X_i]:=\max_{P\in core(\mathbb{V})}E_P[X_i],
$$
for $X_i$ is $\mathcal{F}_i$-measurable, where $core(\mathbb{V}):=\{P \  {\rm additive} | P(A) \leq \mathbb{V}(A), A\in\mathcal{F}\}$.

On one hand, If $\alpha>x_n$, then $\{X_i\geq\alpha\}=\emptyset$, and then $P'_i(X_i\geq\alpha)=\mathbb{V}(X_i\geq\alpha)=0$;
If $x_{k-1}<\alpha\leq x_k$, $k=2,3,\cdots,n$, then $\{X_i\geq\alpha\}=\bigcup\limits_{j\geq k}\{X_i=x_j\}$, and then
$$
\mathbb{V}(X_i\geq\alpha)=\max_{1\leq l\leq m}\sum_{j=k}^{n}p_{lj}=P'_i(X_i\geq\alpha);
$$
If $\alpha<x_1$, then $\{X_i\geq\alpha\}=\Omega$, and then $P'_i(X_i\geq\alpha)=\mathbb{V}(X_i\geq\alpha)=1$, hence for all $\alpha\in\mathbb{R}$, we have
$$
\mathbb{V}(X_i\geq\alpha)=P'_i(X_i\geq\alpha).
$$

On the other hand,
\begin{align*}
P'_i(X_i=x_k)&=\max\limits_{1\leq l\leq m}\sum_{j=k}^{n}p_{lj} - \max\limits_{1\leq l\leq m}\sum_{j=k+1}^{n}p_{lj} \\
&\leq \max\limits_{1\leq l\leq m}p_{lk}=\max\limits_{1\leq l\leq m}P^{l}(X_i=x_k) \\
&=\mathbb{V}(X_i=x_k),
\end{align*}
hence, $P'_i\in core(\mathbb{V})$.

Similarly, given $\varphi_i$, there is a probability $P'_{\varphi_i}$ such that,
$$
\mathbb{V}(\varphi_i(X_i)\geq\alpha)=P'_{\varphi_i}(\varphi_i(X_i)\geq\alpha),
$$
and $P'_{\varphi_i}\in core(\mathbb{V})$.

$\forall n'>0$, Let %$Y_i(\omega):=X_i(\omega_{i})$, $\omega\in\Omega:=\prod_{i=1}^{n'}\Omega_{i}$, $\mathcal{F}:=\prod_{i=1}^{n'}\mathcal{F}_i$, $\bigotimes_{i=1}^{n'}\mathcal{P}_i:=\{\prod_{i=1}^{n'}Q_i|Q_i\in\mathcal{P}_i\}$,
$\left(\Omega, \mathcal{F}, \bigotimes\limits_{i=1}^{n'}\mathcal{P}_i\right):=\left(\prod\limits_{i=1}^{n'}\Omega_{i}, \prod\limits_{i=1}^{n'}\mathcal{F}_i, \left\{\prod\limits_{i=1}^{n'}Q_i|Q_i\in\mathcal{P}_i\right\}\right)$, we have
\begin{align*}
&\mathbb{V}\Big(\sum_{i=1}^{n'}\varphi_i(X_i)\geq \alpha\Big) = \max_{P\in \bigotimes_{i=1}^{n'}\mathcal{P}_i}E_P \Big(I\Big[\sum_{i=1}^{n'}\varphi_i(X_i)\geq \alpha\Big]\Big) \\
&=\max_{P\in \bigotimes_{i=1}^{n'-1}\mathcal{P}_i} \max_{Q\in \mathcal{P}_{n'}} E_P \Big[E_Q \Big[I\Big[\sum_{i=1}^{n'-1}\varphi_i(x_i)+\varphi_{n'}(X_{n'})\geq \alpha\Big]\Big]\Big|_{x_1=X_1,\cdots,x_{n'-1}=X_{n'-1}}\Big] \\
&\leq \max_{P\in \bigotimes_{i=1}^{n'-1}\mathcal{P}_i}  E_P \Big[ \max_{Q\in \mathcal{P}_{n'}} E_Q \Big[I\Big[\sum_{i=1}^{n'-1}\varphi_i(x_i)+\varphi_{n'}(X_{n'})\geq \alpha\Big]\Big]\Big|_{x_1=X_1,\cdots,x_{n'-1}=X_{n'-1}}\Big] \\
&= \mathbb{E} \Big[ \mathbb{E} \Big[I\Big[\sum_{i=1}^{n'-1}\varphi_i(x_i)+\varphi_{n'}(X_{n'})\geq \alpha\Big]\Big]\Big|_{x_1=X_1,\cdots,x_{n'-1}=X_{n'-1}}\Big]\\
&=\mathbb{E}\Big[ \mathbb{V}\Big(\sum_{i=1}^{n'-1}\varphi_i(x_i)+\varphi_{n'}(X_{n'})\geq \alpha\Big)\Big|_{x_1=X_1,\cdots,x_{n'-1}=X_{n'-1}}\Big],
\end{align*}

On the other hand,
\begin{align*}
&\ \ \ \ \mathbb{E}\Big[ \mathbb{V}\Big( \sum_{i=1}^{n'-1}\varphi_i(x_i)+\varphi_{n'}(X_{n'})\geq \alpha\Big)\Big|_{x_1=X_1,\cdots,x_{n'-1}=X_{n'-1}}  \Big] \\
&=\mathbb{E}\Big[ P'_{\varphi_{n'}}\Big( \sum_{i=1}^{n'-1}\varphi_i(x_i)+\varphi_{n'}(X_{n'})\geq \alpha\Big)\Big|_{x_1=X_1,\cdots,x_{n'-1}=X_{n'-1}}  \Big] \\
&=  \max_{P\in \bigotimes_{i=1}^{n'-1}\mathcal{P}_i}  E_P\Big[ P'_{\varphi_{n'}}\Big( \sum_{i=1}^{n'-1}\varphi_i(x_i)+\varphi_{n'}(X_{n'})\geq \alpha\Big)\Big|_{x_1=X_1,\cdots,x_{n'-1}=X_{n'-1}}  \Big]\\
&\leq  \max_{P\in \bigotimes_{i=1}^{n'}\mathcal{P}_i} E_P\Big[ P\Big( \sum_{i=1}^{n'-1}\varphi_i(x_i)+\varphi_{n'}(X_{n'})\geq \alpha\Big)\Big|_{x_1=X_1,\cdots,x_{n'-1}=X_{n'-1}}  \Big]\\
&=  \max_{P\in \bigotimes_{i=1}^{n'}\mathcal{P}_i}P\Big(\sum_{i=1}^{n'}\varphi_i(X_i)\geq \alpha\Big)
\leq \mathbb{V}\Big(\sum_{i=1}^{n'}\varphi_i(X_i)\geq \alpha\Big),
\end{align*}
because of $P'_{\varphi_{n'}}\leq \mathbb{V}$.

Thus, we have for all $\alpha$
$$
\mathbb{V}\Big(\sum_{i=1}^{n'}\varphi_i(X_i)\geq \alpha\Big)=\mathbb{E}\Big[ \mathbb{V}\Big( \sum_{i=1}^{n'-1}\varphi_i(x_i)+\varphi_{n'}(X_{n'})\geq \alpha\Big)\Big|_{x_1=X_1,\cdots,x_{n'-1}=X_{n'-1}}\Big],
$$
i.e. the Fubini independence holds.
\end{example}

\begin{remark}
	The constructions of capacities and upper envelops in Example \ref{ellsberg'smodel} are general for all finite sample space, because $p_{lj}$ can be arbitrary. By the above example, we can construct the upper envelop by core to satisfy the Fubini's independence, in fact, sub-linear expectations can be expressed by upper envelop over its core, see the proof of Theorem 2.1 in Peng \cite{Peng}.
\end{remark}

\subsection{Infinite sample case}
Next example shows that the Ellsberg's model (with infinite sample space) also satisfies the Fubini independence for $\mathcal{U}_{\hat{f}}$ under upper envelop $\mathbb{E}$.
\begin{example}
As considered in example \ref{ellsberg'smodel}, for every $i\geq 1$, $\Omega_i$ is stand for the $i$-th urn.

Step 1:

Let $\{A_k^i, k=1,\ldots,n\}$ are the finite divisions of $\Omega_i$. $\forall \omega \in \Omega_i$, $X_i(\omega):=\sum\limits_{k=1}^{n}x_kI_{A_k^i}(\omega)$, and $x_1\leq x_2\leq\cdots\leq x_n$, if not we just reorder divisions of $\Omega_i$.
Denote $P^l(X_i=x_j):=p_{lj}$, for $l=1,2,\cdots,m$, $j=1,2,\cdots,n$, where $p_{lj}\geq 0$, $\sum\limits_{j=1}^{n}p_{lj}=1$, and $\mathcal{P}_i:=\{P^l | l=1,2,\cdots,m\}$. Then we can get the conclusion by the example \ref{ellsberg'smodel}, using the same method.

Step 2:

Let $X_i$ is an arbitrary r.v. in $(\Omega_i, \mathcal{F}_i)$, we can find a sequence $\{X_i^r, r\geq 1\}$ of increasing step functions in Step 1 to approximate $X_i$.

Given bounded continuous $\varphi_i$, we can get
$$
\mathbb{V}(\varphi_i(X_i)\geq\alpha)=\lim_{r\rightarrow \infty}\mathbb{V}(\varphi_i(X_i^r)\geq\alpha),
$$ by Dominated Convergence Theorem (Thm 8.9 in Denneberg), and
$$
P'_{\varphi_i}(\varphi_i(X_i)\geq\alpha)=\lim_{r\rightarrow \infty}P'_{\varphi_i}(\varphi_i(X_i^r)\geq\alpha),
$$
by classic monotonic convergence theorem.

Therefore $
\mathbb{V}(\varphi_i(X_i)\geq\alpha)=P'_{\varphi_i}(\varphi_i(X_i)\geq\alpha)
$ by Step 1. And we can get the conclusion.
\end{example}

\section{Weak Law of Large Numbers under Exponential Independence}

In this section, we let $\mathbb{E}[X]:=\sup\limits_{P\in\mathcal{P}}E_P[X]$, $\mathcal{E}[X]:=\inf\limits_{P\in\mathcal{P}}E_P[X]$, $\mathbb{V}(A):=\sup\limits_{P\in\mathcal{P}}P(A)$, $\nu(A):=\inf\limits_{P\in\mathcal{P}}P(A)$, where $\mathcal{P}$ is a set of probabilities.

Strong law of large numbers under Exponential independence can see Chen, Huang, Wu \cite{CHW}.

\begin{thm}\label{ThmWLLN}
	Let $\{X_n\}_{n\geq 1}$ be a sequence of exponential independent random variables, and $\mathbb{E}[X_k]=\overline{\mu}$, $\mathcal{E}[X_k]=\underline{\mu}$, $\forall k\geq 1$, $\sup\limits_{k\geq 1}\mathbb{E}[|X_k|^{1+\alpha}]<\infty$, for some $0<\alpha<1$. Then we have
	\begin{equation}\label{WLLN}
	\lim_{n\rightarrow\infty}\nu\left(\underline{\mu}-\epsilon\leq \frac{1}{n}\sum_{k=1}^n X_k \leq\overline{\mu}+ \epsilon\right)=1, \ \ \forall \epsilon>0.
	\end{equation}
\end{thm}
\Proof. If
\begin{equation}\label{WLLN1}
\lim_{n\rightarrow\infty}\nu\left(\frac{1}{n}\sum_{k=1}^n X_k \leq\overline{\mu}+ \epsilon\right)=1
\end{equation}
holds, then we consider $-X_k$ instead of $X_k$,
\begin{equation*}
\lim_{n\rightarrow\infty}\nu\left(\frac{1}{n}\sum_{k=1}^n (-X_k) \leq\mathbb{E}[-X_k]+ \epsilon\right)=1,
\end{equation*}
i.e.
\begin{equation}\label{WLLN2}
\lim_{n\rightarrow\infty}\nu\left(\frac{1}{n}\sum_{k=1}^n X_k \geq\underline{\mu}- \epsilon\right)=1,
\end{equation}
combined (\ref{WLLN1}) and (\ref{WLLN2}) we have (\ref{WLLN}), in fact
\begin{align*}
&\ \ \ \ \nu\left(\underline{\mu}-\epsilon\leq \frac{1}{n}\sum_{k=1}^n X_k \leq\overline{\mu}+ \epsilon\right) \\
&=1-\mathbb{V}\left( \left[\frac{1}{n}\sum_{k=1}^n X_k-\overline{\mu}>\epsilon\right]\cup \left[\frac{1}{n}\sum_{k=1}^n X_k-\underline{\mu}<-\epsilon\right] \right)\\
&\geq 1-\left( \mathbb{V}\left( \frac{1}{n}\sum_{k=1}^n X_k-\overline{\mu}>\epsilon\right) + \mathbb{V}\left(\frac{1}{n}\sum_{k=1}^n X_k-\underline{\mu}<-\epsilon\right) \right)\\
&=\nu\left(\frac{1}{n}\sum_{k=1}^n X_k \leq\overline{\mu}+ \epsilon\right)+ \nu\left(\frac{1}{n}\sum_{k=1}^n X_k \geq\underline{\mu}- \epsilon\right)-1\\
&\rightarrow 1+1-1=1.
\end{align*}
Hence we only need to prove (\ref{WLLN1}).

Let $S_n:=\sum\limits_{k=1}^{n}X_k$, we defined
\begin{equation*}
f_n(x):=\left(-\frac{n}{\log(1+n)}\right)\vee\left( x\wedge \frac{n}{\log(1+n)} \right),
\end{equation*}
\begin{equation*}
\hat{f_n}(x):=x-f_n(x),
\end{equation*}
\begin{equation*}
\overline{X}_n:=f_n(X_n-\overline{\mu})-\mathbb{E}[f_n(X_n-\overline{\mu})]+\overline{\mu}.
\end{equation*}
Then $\mathbb{E}[\overline{X}_n]=\overline{\mu}$, $\forall n\geq 1$. Let $\overline{S}_n:=\sum\limits_{k=1}^n \overline{X}_k$, we have
\begin{equation}\label{jieduan}
\frac{S_n-n\overline{\mu}}{n}=\frac{S_n-\overline{S}_n}{n}+\frac{\overline{S}_n-n\overline{\mu}}{n}.
\end{equation}

We first prove that, $\forall \epsilon>0$ we have
\begin{equation}\label{WLLN3}
\lim_{n\rightarrow\infty}\mathbb{V}\left(\frac{\overline{S}_n-n\overline{\mu}}{n}>\epsilon\right)=0,
\end{equation}
by Chebyshev's inequality we have
\begin{align*}
    \mathbb{V}\left(\frac{\overline{S}_n-n\overline{\mu}}{n}>\epsilon\right)&=\mathbb{V}\left(\sum_{k=1}^n\left(\overline{X}_k-\mathbb{E}[\overline{X}_k]\right)>n\cdot\epsilon\right) \\
&=\mathbb{V}\left(\frac{m\log(1+n)}{n}\sum_{k=1}^n\left(\overline{X}_k-\mathbb{E}[\overline{X}_k]\right)>m\log(1+n)\cdot\epsilon\right) \\
&\leq \frac{\mathbb{E}\left[e^{\frac{m\log(1+n)}{n}\sum_{k=1}^n\left(\overline{X}_k-\mathbb{E}[\overline{X}_k]\right)}\right]}{e^{m\epsilon\log(1+n)}}\\
&\leq \frac{\sup_{n\geq 1} \mathbb{E}\left[e^{\frac{m\log(1+n)}{n}\sum_{k=1}^n\left(\overline{X}_k-\mathbb{E}[\overline{X}_k]\right)}\right]}{e^{m\epsilon\log(1+n)}}.
\end{align*}
By Lemma 3.1 in \cite{CWL} \footnote{The independence condition in Lemma 3.1 in \cite{CWL} is not the exponential independence, but it is easy to prove it under exponential independence without change the proof.} we have, $\sup\limits_{n\geq 1} \mathbb{E}\left[e^{\frac{m\log(1+n)}{n}\sum_{k=1}^n\left(\overline{X}_k-\mathbb{E}[\overline{X}_k]\right)}\right]<\infty$, $\forall m>1$. Meanwhile, $\lim\limits_{n\rightarrow\infty}e^{m\epsilon\log(1+n)}=\infty$, hence we have (\ref{WLLN3}).

Secondly, we prove $\forall \epsilon>0$,
\begin{equation}\label{WLLN4}
\lim_{n\rightarrow\infty}\mathbb{V}\left(\frac{S_n-\overline{S}_n}{n}>\epsilon\right)=0.
\end{equation}
By the reason of
$$X_n=\overline X_n+\hat{f}_n(X_n-\overline\mu)+\mathbb{E}\left[f_n(X_n-\overline\mu)\right],$$
we have
$$\frac{1}{n}S_n= \frac{1}{n}\overline S_n+ \frac 1n
\sum_{k=1}^n \hat{f_k}(X_k-\overline
\mu)+\frac1n
\sum_{k=1}^n\mathbb{E}\left[f_k(X_k-\overline
\mu)\right].$$
By the subadditivity and translation invariance of $\mathbb{E}[\cdot]$, we obtain
\begin{eqnarray*}
	\mathbb{E}\left[f_k(X_k-\overline
	\mu)\right]&=&\mathbb{E}\left[(X_k-\overline{\mu})-\hat{f_k}(X_k-\overline
	\mu)\right]\\
	&\leq&\mathbb{E}[X_k]-\overline{\mu}+\mathbb{E}\left[-\hat{f_k}(X_k-\overline
	\mu)\right]\\
	&\leq& \mathbb{E}\left[|\hat{f_k}(X_k-\overline
	\mu)|\right].
\end{eqnarray*}
Notice that $$|\hat{f}_k(X_k-\overline{\mu})|=|X_k-\overline{\mu}-f_k(X_k-\overline{\mu})|\leq |X_k-\overline
\mu|I_{\{|X_k-\overline\mu|>\frac{k}{\log(k+1)}\}},$$
therefore,
\begin{equation}\label{R}
\frac{1}{n}S_n\leq \frac{1}{n}\overline S_n+ \frac 1n
\sum_{k=1}^n|X_k-\overline
\mu|I_{\{|X_k-\overline\mu|>\frac{k}{\log(1+k)}\}}+\frac1n
\sum_{k=1}^n\mathbb{E}\left[|X_k-\overline
\mu|I_{\{|X_k-\overline\mu|>\frac{k}{\log(1+k)}\}}\right].
\end{equation}
Applying the H\"{o}lder and Chebyshev inequality, we have
\begin{eqnarray*}
	&&\sum\limits_{k=1}^\infty\frac{\mathbb{E}\left[|X_k-\overline
		\mu|I_{\{|X_k-\overline\mu|>\frac{k}{\log(1+k)}\}}\right]}{k}\\
	&\leq& \sum\limits_{k=1}^\infty\frac 1k(\mathbb{E}[|X_k-\overline{\mu}|^{1+\alpha}])^{1/(1+\alpha)}(\mathbb{E}[I_{\{|X_k-\overline\mu|>\frac{k}{\log(1+k)}\}}])^{\alpha/(1+\alpha)}\\
	&\leq& \sum\limits_{k=1}^\infty\frac{[\log(1+k)]^\alpha}{
		k^{1+\alpha }} \mathbb{E}[|X_k-\overline \mu|^{1+\alpha}]\\
	&\leq&\sup\limits_{k\geq1}\mathbb{E}[|X_k-\overline \mu|^{1+\alpha}]\sum\limits_{k=1}^\infty\frac{[\log(1+k)]^\alpha}
	{k^{1+ \alpha }}<\infty.\end{eqnarray*}
By Kronecker Lemma, we have
\begin{equation}\label{E}
\lim_{n\to\infty}\frac 1n \sum_{k=1}^n\mathbb{E}\left[|X_k-\overline
\mu|I_{\{|X_k-\overline\mu|>\frac{k}{\log(1+k)}\}}\right]=0.
\end{equation}
Hence by (\ref{R}) and (\ref{E}) we have
\begin{align*}
&\ \ \ \ \mathbb{V}\left(\frac{S_n-\overline{S}_n}{n}>\epsilon\right) \\
&\leq \mathbb{V}\left(\frac 1n \sum_{k=1}^n\left(|X_k-\overline \mu|I_{\{|X_k-\overline\mu|>\frac{k}{\log(1+k)}\}}+\mathbb{E}\left[|X_k-\overline \mu|I_{\{|X_k-\overline\mu|>\frac{k}{\log(1+k)}\}}\right]\right)>\epsilon\right) \\
&\leq \frac{\sum_{k=1}^n\mathbb{E}\left[ |X_k-\overline \mu|I_{\{|X_k-\overline\mu|>\frac{k}{\log(1+k)}\}}+\mathbb{E}\left[|X_k-\overline \mu|I_{\{|X_k-\overline\mu|>\frac{k}{\log(1+k)}\}}\right] \right]}{n\epsilon} \\
&=\frac{2\sum_{k=1}^n\mathbb{E}\left[ |X_k-\overline \mu|I_{\{|X_k-\overline\mu|>\frac{k}{\log(1+k)}\}} \right]}{n\epsilon}\rightarrow 0,
\end{align*}
i.e. we get (\ref{WLLN4}).

By (\ref{jieduan}), (\ref{WLLN3}) and (\ref{WLLN4}), we have
\begin{align*}
&\ \ \ \ \nu\left(\frac{1}{n}\sum_{k=1}^n X_k-\overline{\mu}\leq \epsilon\right) \\
&=\nu\left(\frac{S_n-n\overline{\mu}}{n}\leq \epsilon\right)\\
&=1-\mathbb{V}\left(\frac{S_n-n\overline{\mu}}{n}> \epsilon\right) \\
&\geq 1-\left(\mathbb{V}\left(\frac{S_n-\overline{S}_n}{n}> \frac{\epsilon}{2}\right) + \mathbb{V}\left(\frac{\overline{S}_n-n\overline{\mu}}{n}> \frac{\epsilon}{2}\right)\right) \\
&\rightarrow 1-(0+0)=1.
\end{align*}
Then we have (\ref{WLLN1}), and then (\ref{WLLN}).

\section{Simulations}

In this section, two different Ellsberg's models are considered to demonstrate how the weak law of large numbers under Exponential Independence works.

\subsection{$X_i$ only has mean uncertainty}

\begin{figure}[H]
\centering  %图片全局居中
\subfigure{
\label{Fig.sub.1.1}
\includegraphics[width=0.45\textwidth]{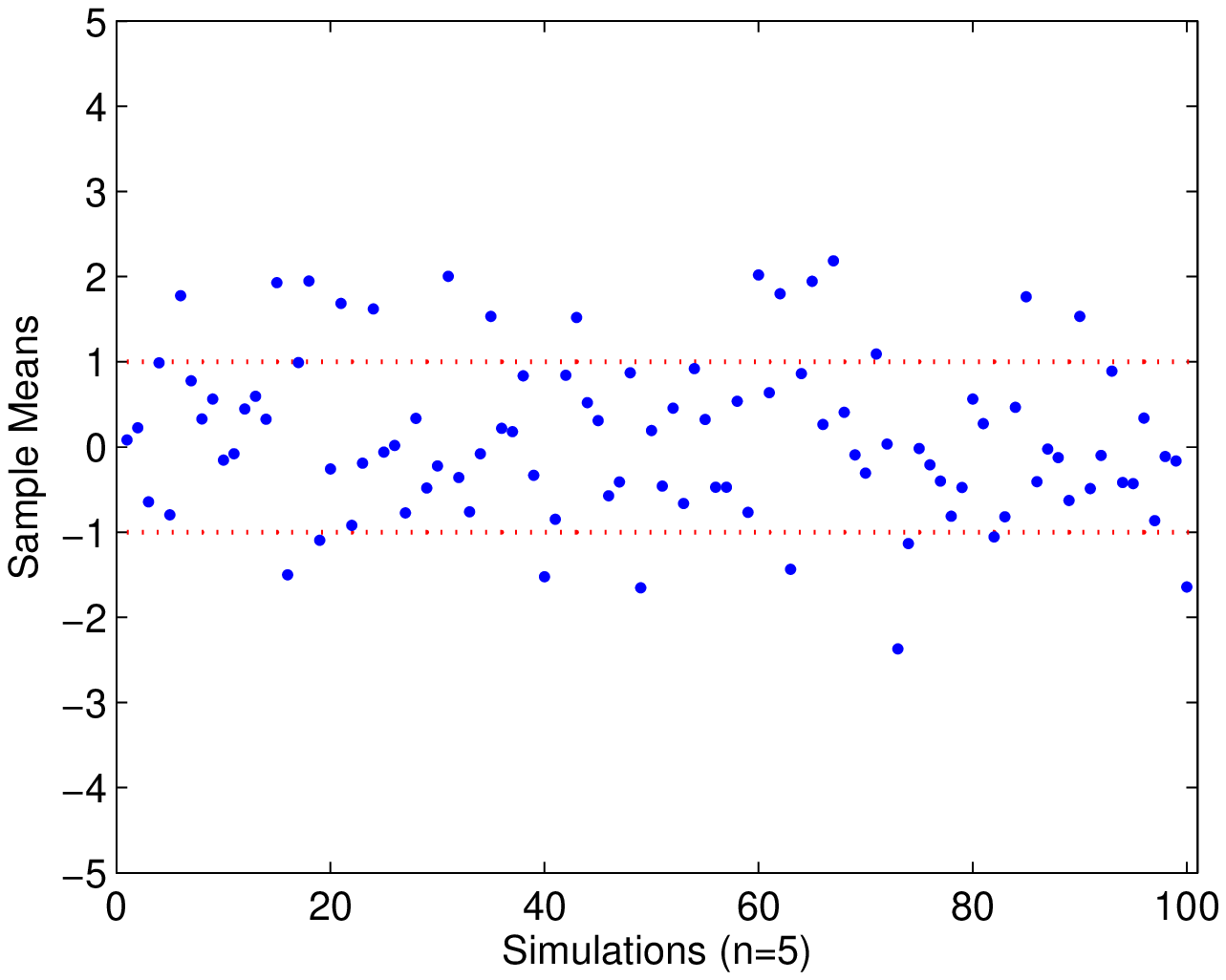}}
\subfigure{
\label{Fig.sub.1.2}
\includegraphics[width=0.45\textwidth]{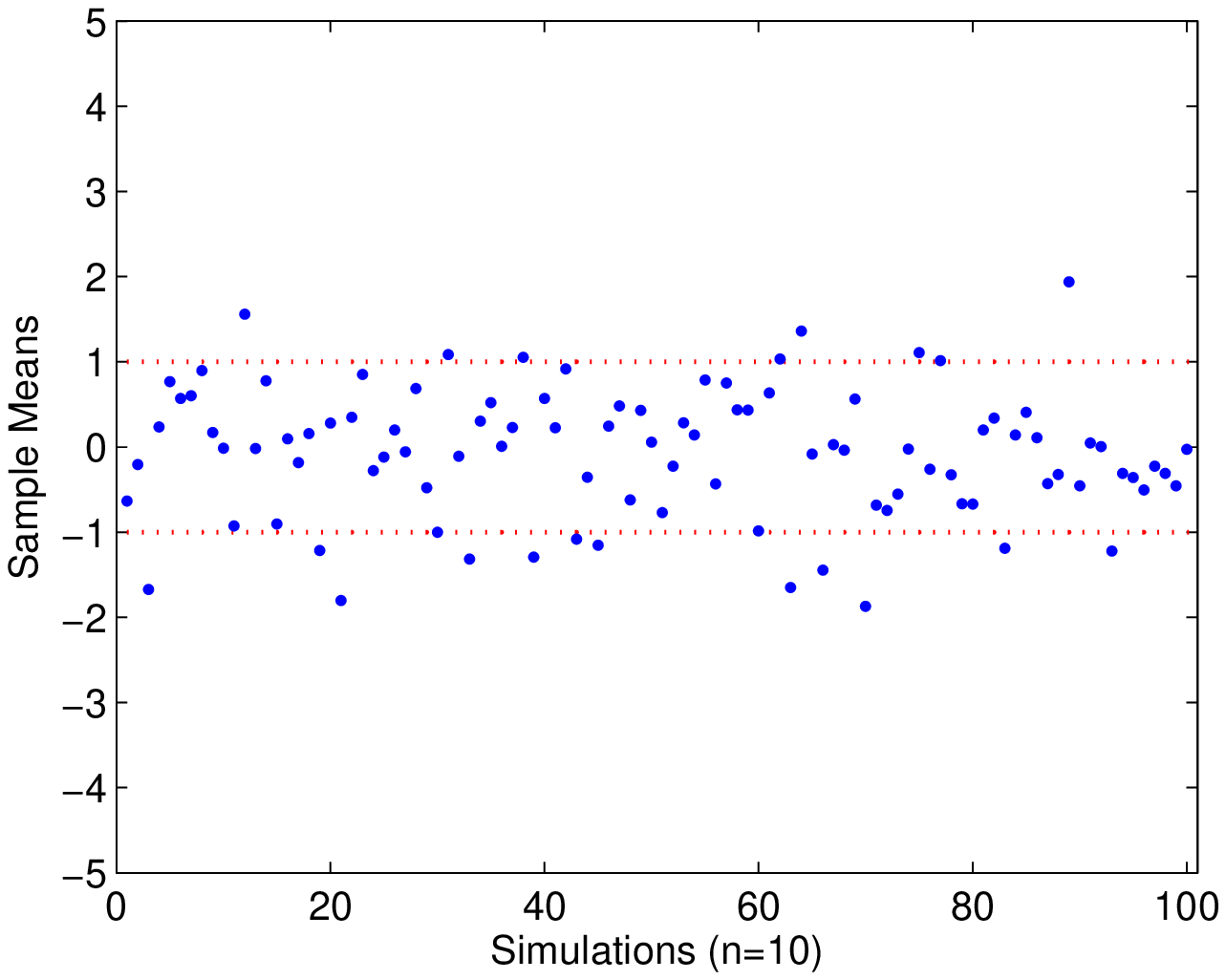}}
\subfigure{
\label{Fig.sub.1.3}
\includegraphics[width=0.45\textwidth]{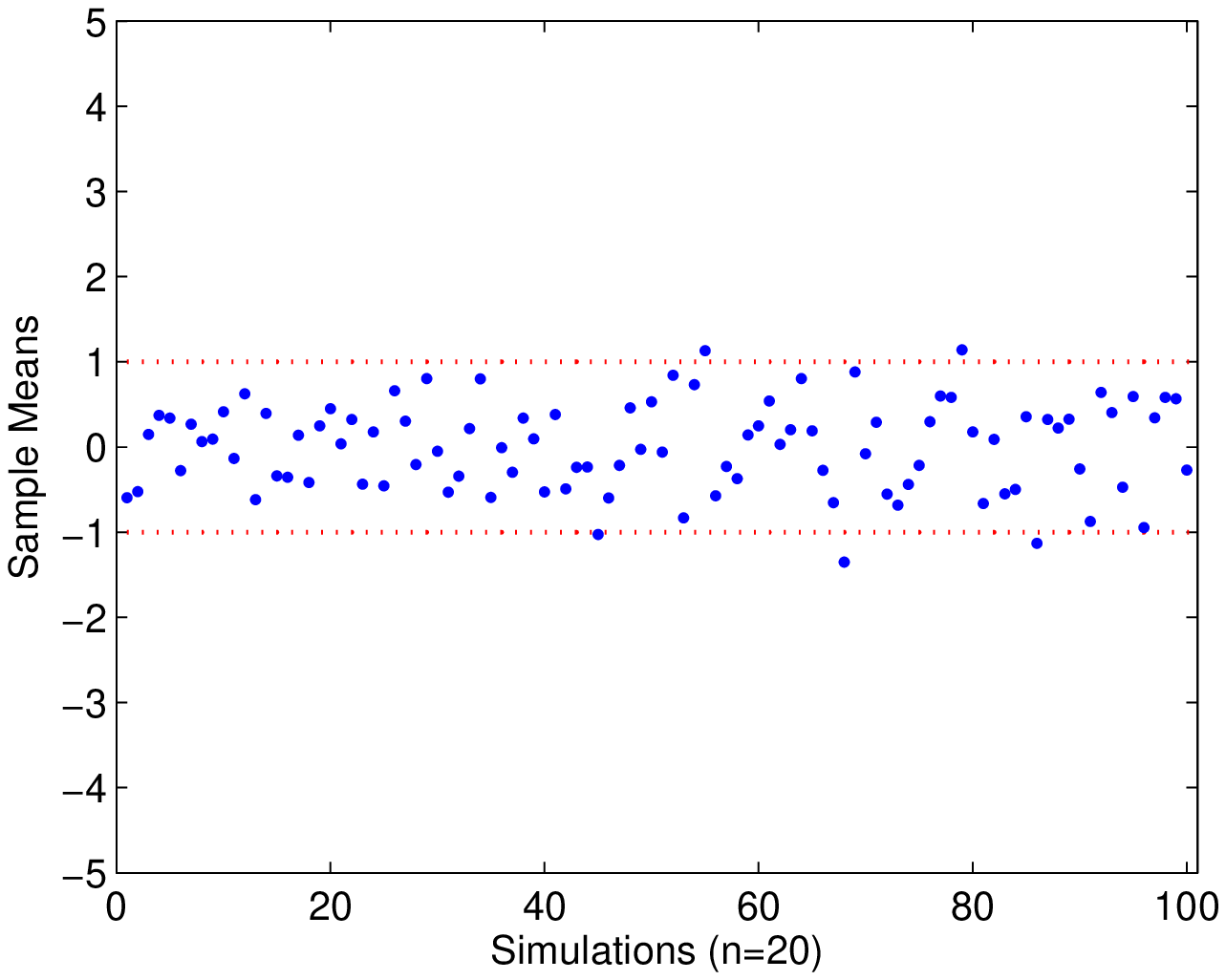}}
\subfigure{
\label{Fig.sub.1.4}
\includegraphics[width=0.45\textwidth]{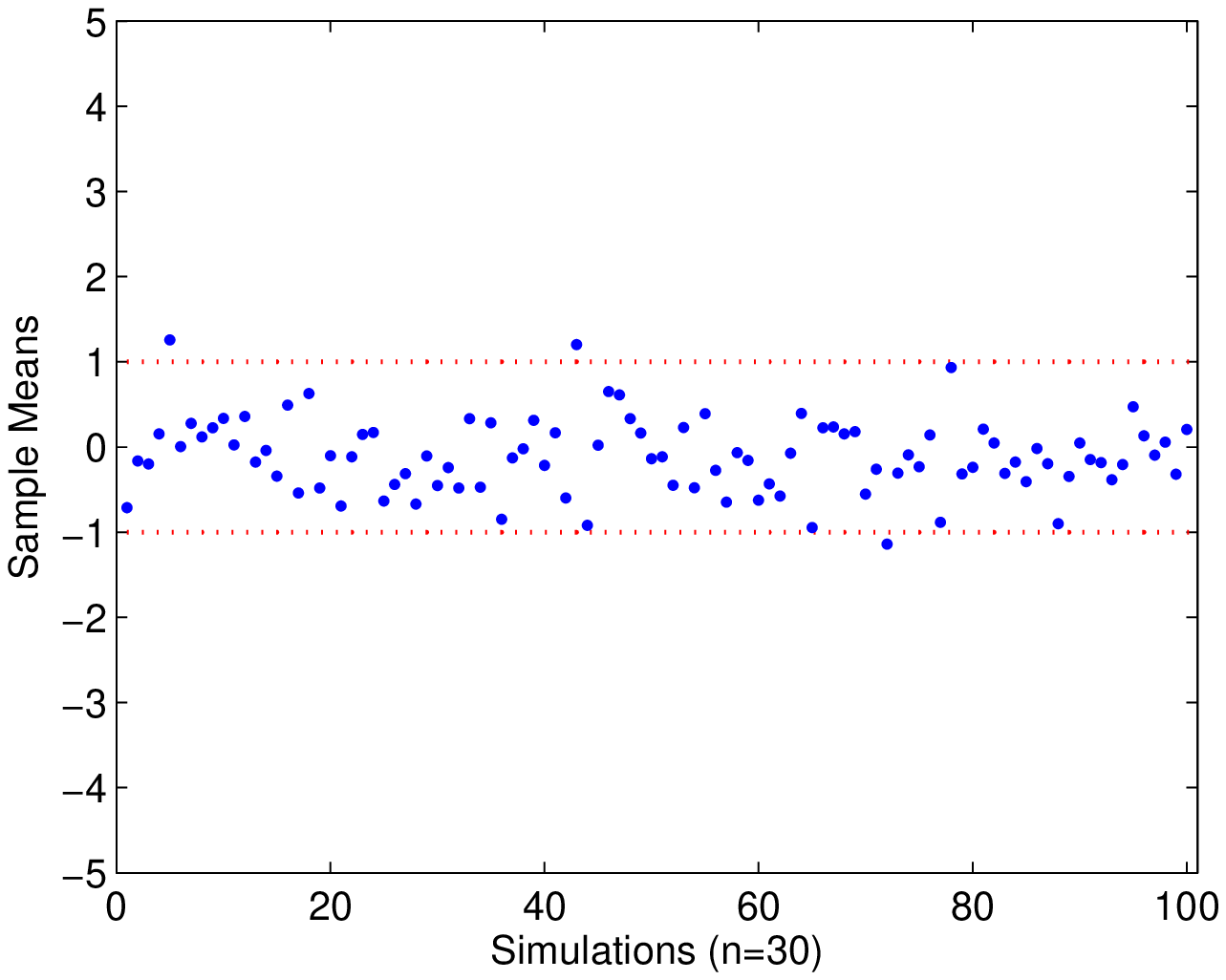}}
\subfigure{
\label{Fig.sub.1.5}
\includegraphics[width=0.45\textwidth]{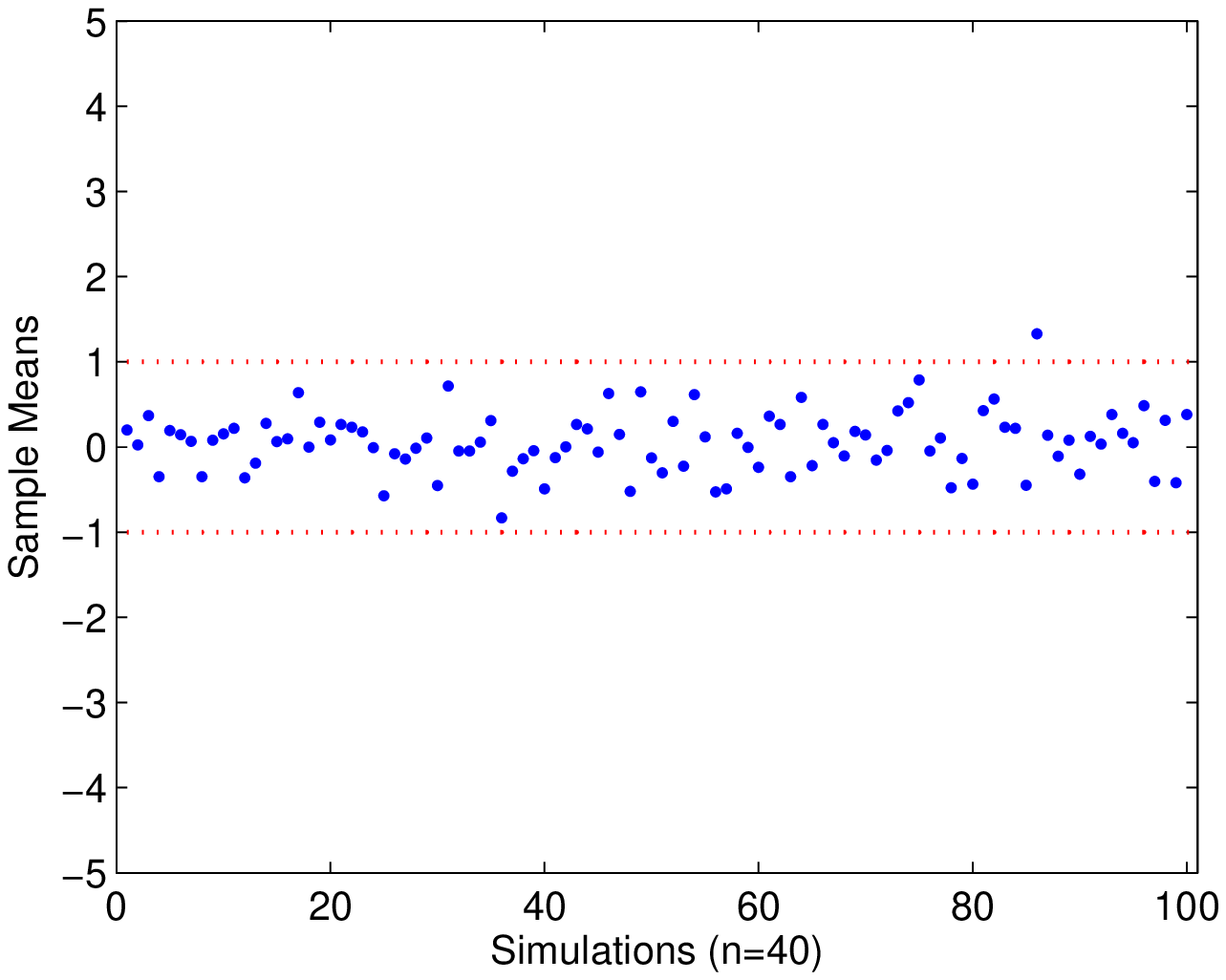}}
\subfigure{
\label{Fig.sub.1.6}
\includegraphics[width=0.45\textwidth]{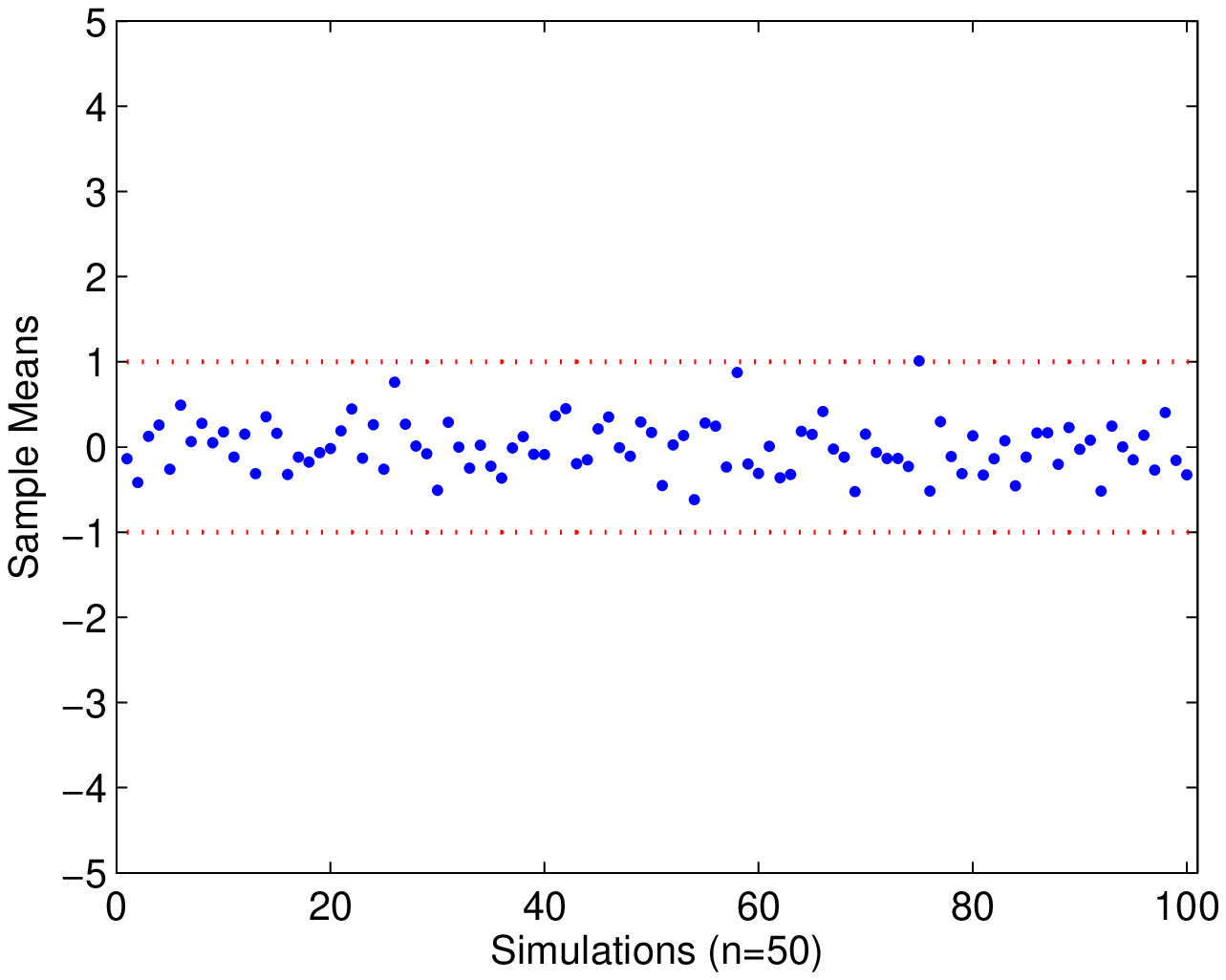}}
\caption{Sample means in 100 simulations}
\label{Fig.main1}
\end{figure}

First, we consider a sequence of Ellsberg's urns $\{X_i\}_{i=1}^n$ which satisfy normal distributions with  $-1\leq\mathbb{E}[X_i]\leq1$ and determined standard deviation $\sigma=2$. Set sample size $n=5, 10,  20, 30, 40, 50$. In every scenario, we repeat simulation 100 times and obtain the sample mean in every times respectively. The six subfigures in Fig \ref{Fig.main1} show that with the growth of sample size, an increasing number of sample means, represented by blue points in each subfigure, lie between the lower mean and the upper mean.

%\begin{figure}[H]
%\centering  %图片全局居中
%\subfigure[name1]{
%\label{Fig.sub.1}
%\includegraphics[width=0.45\textwidth]{n=5.eps}}
%\subfigure[name2]{
%\label{Fig.sub.2}
%\includegraphics[width=0.45\textwidth]{n=10.eps}}
%\subfigure[name3]{
%\label{Fig.sub.3}
%\includegraphics[width=0.45\textwidth]{n=20.eps}}
%\subfigure[name4]{
%\label{Fig.sub.4}
%\includegraphics[width=0.45\textwidth]{n=30.eps}}
%\subfigure[name5]{
%\label{Fig.sub.5}
%\includegraphics[width=0.45\textwidth]{n=40.eps}}
%\subfigure[name6]{
%\label{Fig.sub.6}
%\includegraphics[width=0.45\textwidth]{n=50.eps}}
%\caption{Main name}
%\label{Fig.main}
%\end{figure}

\begin{figure}[H] %H为当前位置，!htb为忽略美学标准，htbp为浮动图形
\centering %图片居中
\includegraphics[width=0.45\textwidth]{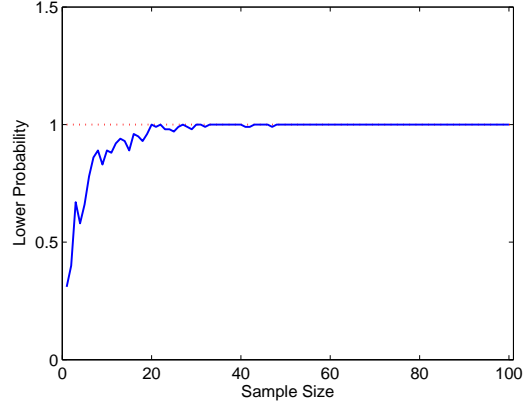} % 插入图片，[]中设置图片大小，{}中是图片文件名
\caption{Lower probability in different sample size} %最终文档中希望显示的图片标题
\label{Fig.main2} %用于文内引用的标签
\end{figure}

Then, the lower probability $v$ that sample means lie between the lower mean and the upper mean can be calculated as follows. We claim that every $P\in\mathcal{P}$ are equivalent in this model. Indeed, every Ellsberg's urns $\{X_i\}_{i=1}^n$ can be characterised by the following stochastic differential equations,
$$
dX_i(t)=\mu_{i}(t)dt+\sigma_{i}dB_{t}^{P},
$$
where $P\in\mathcal{P}$, $B_{t}^{P}$ is a standard Brownian Motion under $P$, $\mu_i(t)$ is an adapted process and $\sigma_i\equiv 2$ for $i=1,\ldots,n$. For any fixed $Q\in\mathcal{P}$, by Girsanov transformation, we have
$$
dX_i(t)=(\mu_{i}(t)-\sigma_{i}\theta^{(Q)}(t))dt+\sigma_{i}dB_{t}^{Q},
$$
where $dB_t^{Q}=\theta^{(Q)}(t)dt + dB_t^{P}$, $B_{t}^{Q}$ is a standard Brownian Motion under $Q$ and $\theta^{(Q)}(t)$ is an adapted process.
This proves our claim, and we also have $-1 \leq \mu_{i}(t)-\sigma_{i}\theta^{(Q)}(t)\leq 1$ by the condition that $-1\leq\mathbb{E}[X_i]\leq1$.
Therefore, we can use the frequency at which sample means lie between the lower mean and the upper mean as the approximation of the lower probability $v$. Fig \ref{Fig.main2} shows that the lower probability $v$ that sample means lie between the lower mean and the upper mean tends to be one with $n\rightarrow\infty$, which satisfies (\ref{WLLN}) in Theorem \ref{ThmWLLN}.

\subsection{$X_i$ has both mean uncertainty and variance uncertainty}

\begin{figure}[H]
\centering  %图片全局居中
\subfigure{
\label{Fig.sub.3.1}
\includegraphics[width=0.45\textwidth]{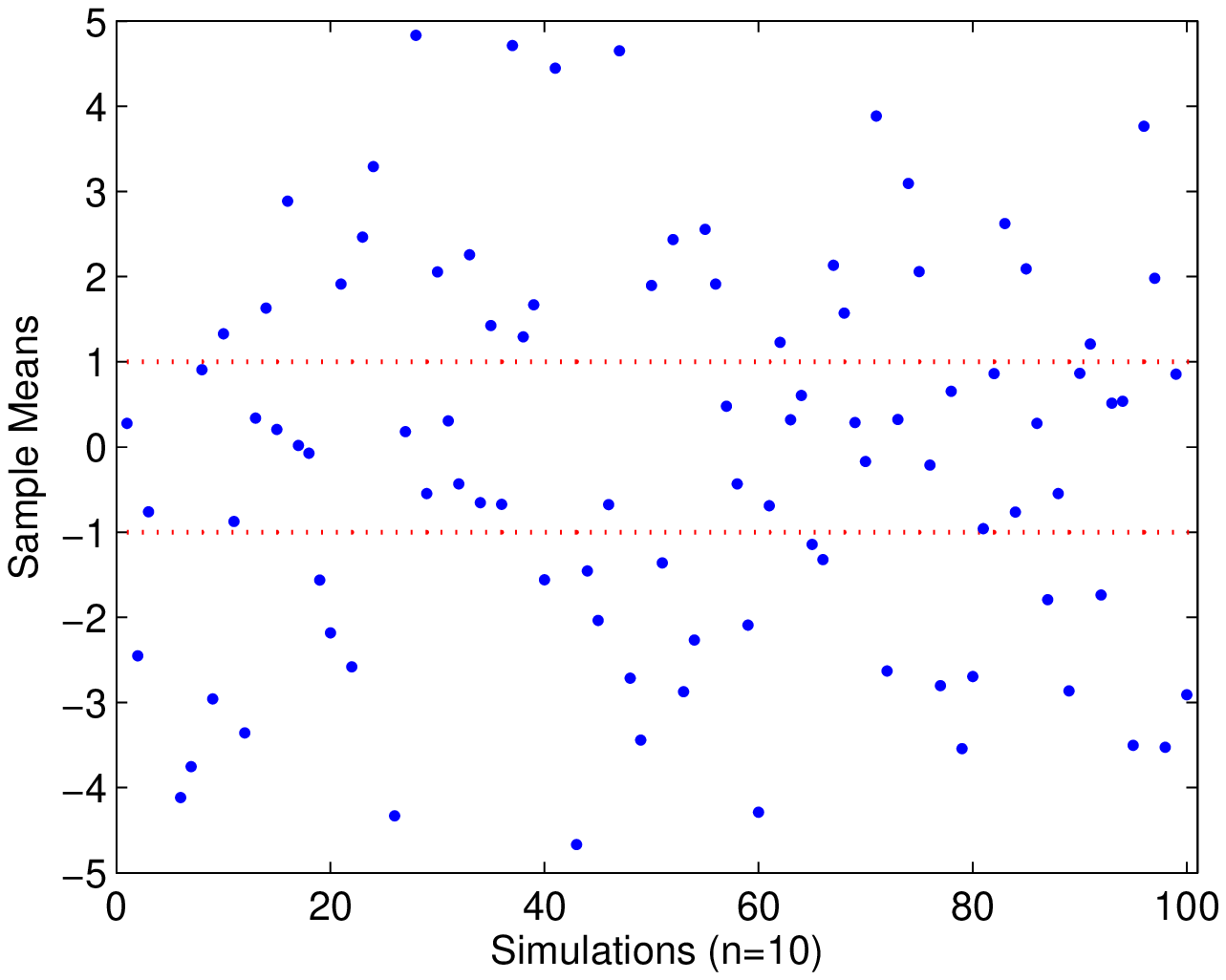}}
\subfigure{
\label{Fig.sub.3.2}
\includegraphics[width=0.45\textwidth]{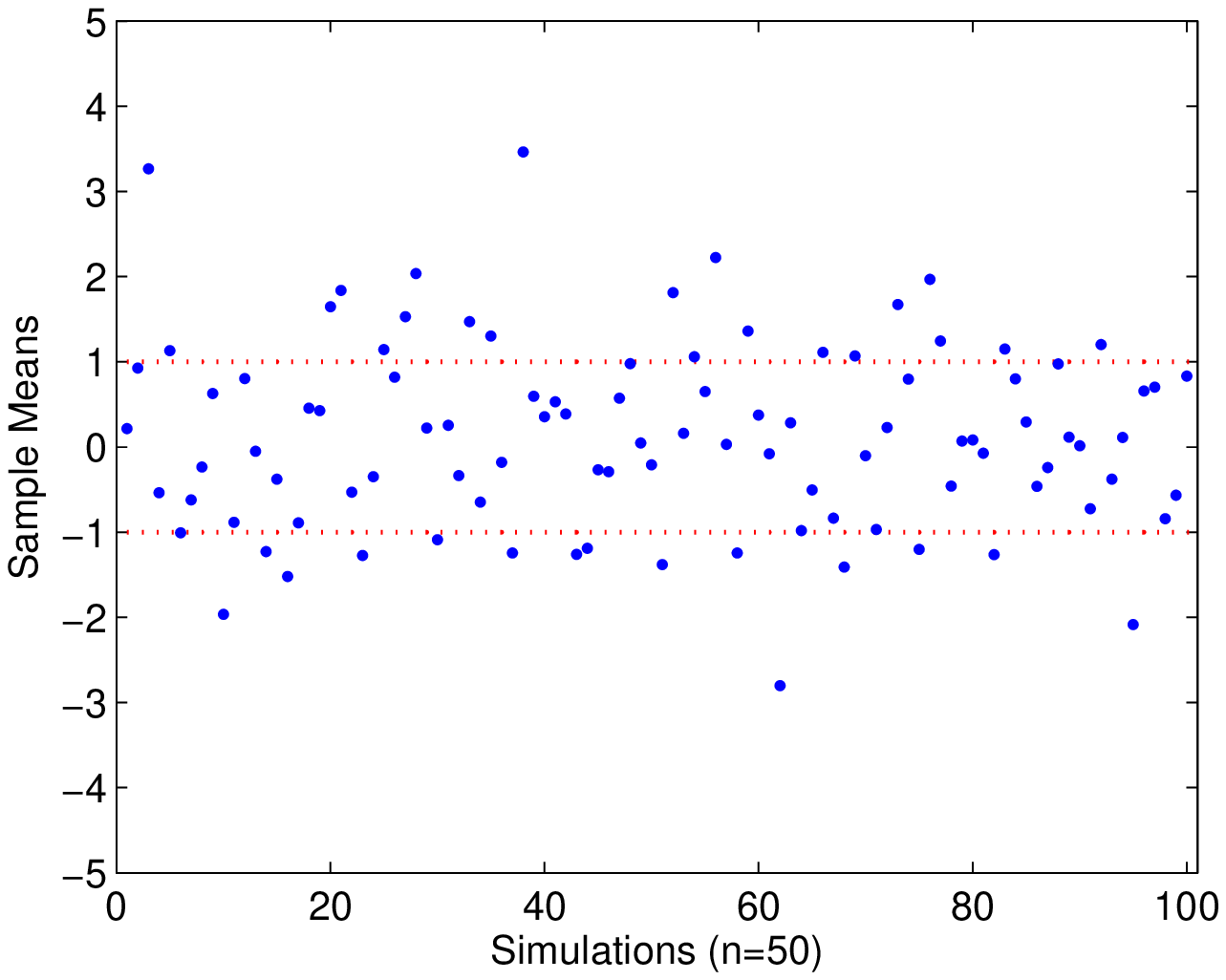}}
\subfigure{
\label{Fig.sub.3.3}
\includegraphics[width=0.45\textwidth]{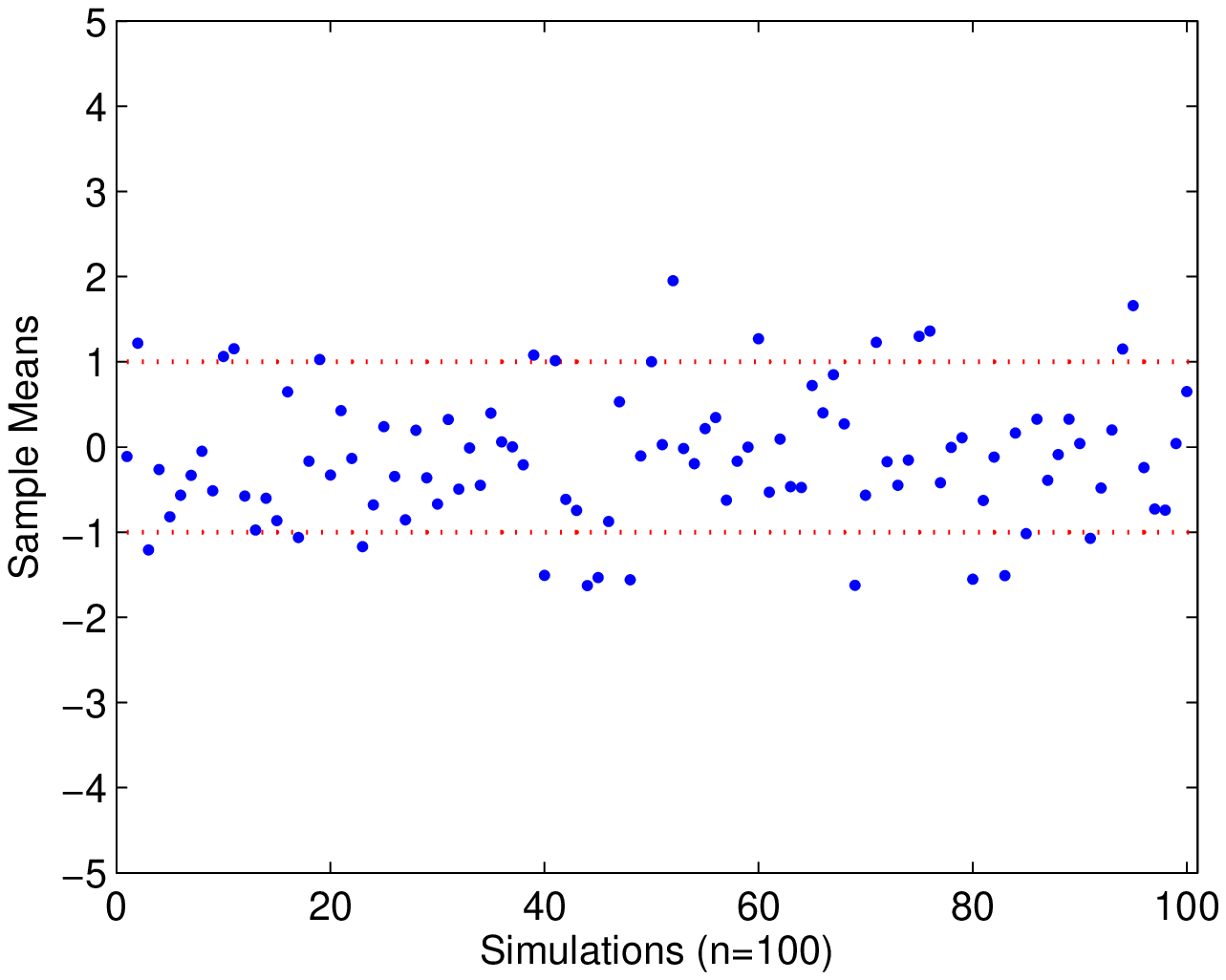}}
\subfigure{
\label{Fig.sub.3.4}
\includegraphics[width=0.45\textwidth]{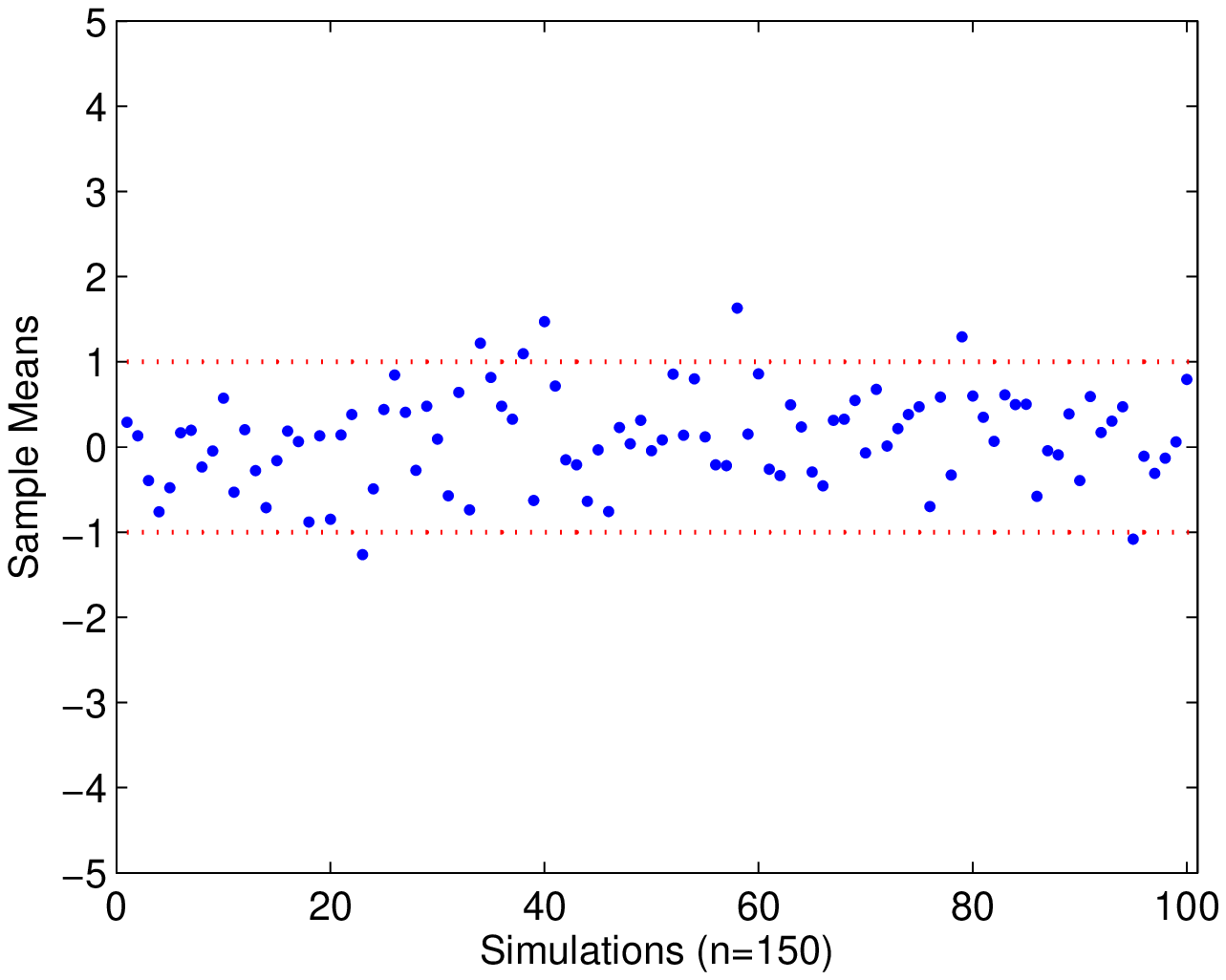}}
\subfigure{
\label{Fig.sub.3.5}
\includegraphics[width=0.45\textwidth]{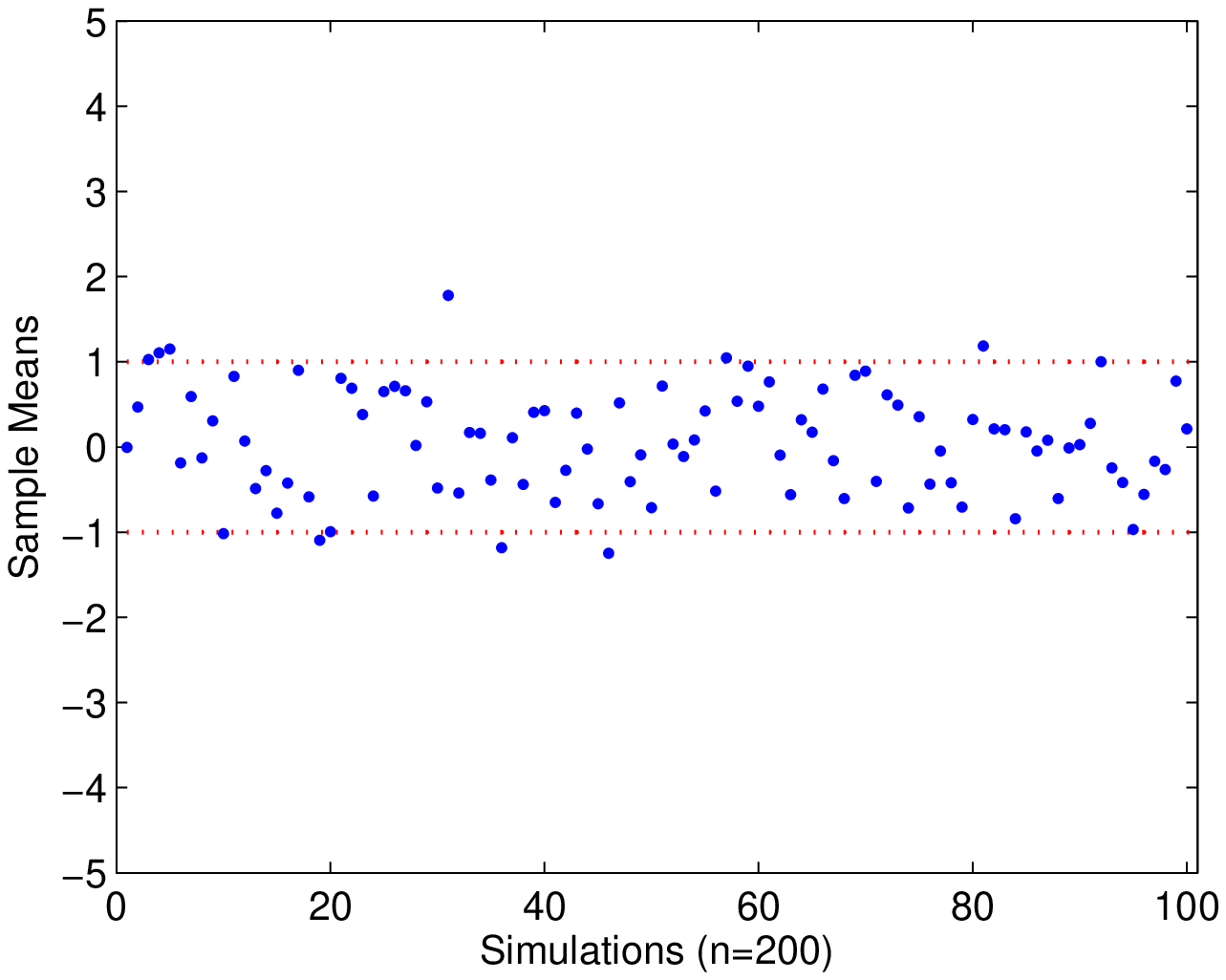}}
\subfigure{
\label{Fig.sub.3.6}
\includegraphics[width=0.45\textwidth]{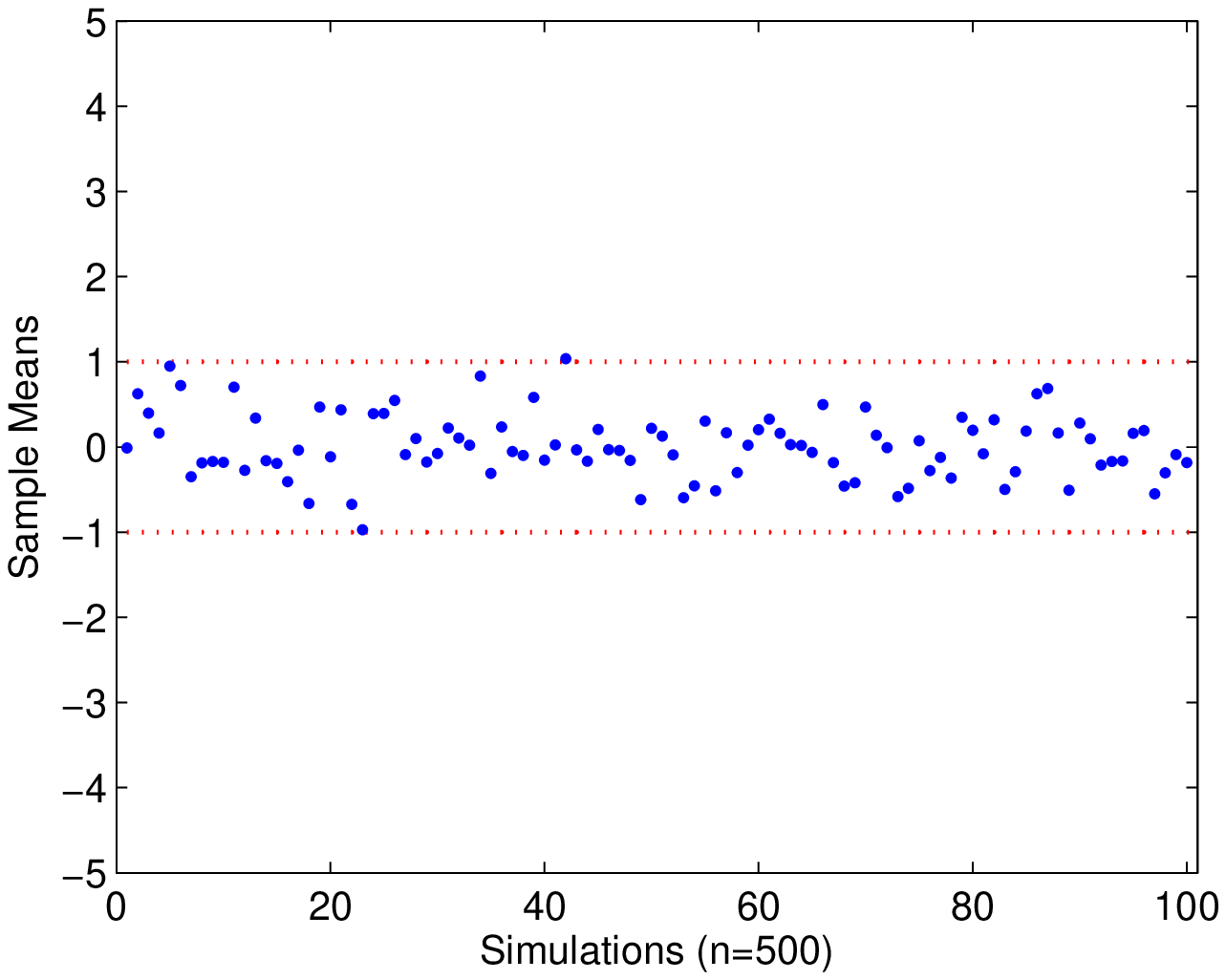}}
\caption{Sample means in 100 simulations}
\label{Fig.main3}
\end{figure}

In this subsection, we consider another sequence of Ellsberg's urns $\{X_i\}_{i=1}^n$ which satisfy normal distributions with  $-1\leq\mathbb{E}[X_i]\leq1$ and uncertain standard deviation $5\leq\sigma_i\leq10$. Set sample size $n=10,  50, 100, 150, 200, 500$. In every scenario, we repeat simulation 100 times and obtain the sample mean in every times respectively. The six subfigures in Fig \ref{Fig.main3} show our simulation results. And we can obtain the similar conclusion as above.

\appendix
\section*{Proof of Theorem \ref{fubiniduli}\label{Appendix}}

\Proof. For $p\in\mathbb{N}$ and $r\in\mathbb{R}_{+}$ define the transformation
$$
u_{p}(r)=\sup\left\{\frac{i}{2^p} : i\in\mathbb{Z}_{+}, \ r\geq\frac{i}{2^p} \right\}.
$$
For $\hat{f}(x_1,x_2)=e^{\varphi_1(x_1)+\varphi_2(x_2)}$, let $f_p(x_1,x_2) = u_p [ \hat{f}(x_1,x_2)]$, as we can find a finite chain of sets $\{ \{ e^{\varphi_1(x_1)+\varphi_2(x_2)}\geq \frac{i}{2^p} \} |i=1,\cdots,n\}:=\{A_1, A_2, \cdots, A_n\}\subset\mathcal{U}_{\hat{f}}$, $n\geq 1$, such that $A_{i+1}\subset A_i$ for $i=1, 2,\cdots, n-1$, $0\leq\alpha_1<\alpha_2<\cdots<\alpha_n$, and $\alpha_i=\frac{i}{2^p}$, %and $\alpha_i=\frac{i}{2^p}$.
\begin{equation}\label{shizi2}
f_p(x_1,x_2)=\alpha_1 I_{A_1}(x_1,x_2) + \sum_{i=2}^{n}(\alpha_{i}-\alpha_{i-1})I_{A_i}(x_1,x_2),
\end{equation}
%$\forall y\in Y$, $g(\cdot, y)$ form a class of comonotonic functions,
where $\alpha_{i}-\alpha_{i-1}>0$, $A_i=\{ (x_1,x_2)|e^{\varphi_1(x_1)+\varphi_2(x_2)}\geq \frac{i}{2^p} \}$. In fact, we can find a finite chain because of the boundness of $\varphi_i$. In this case, $\mathcal{U}_{f_p}=\{A_1, A_2, \cdots, A_n, \mathbb{R}^2, \emptyset\}$. The set $\{f_p, I_{A_1}, I_{A_2}, \cdots, I_{A_n}\}$ forms a comonotonic class. In fact, every pair $I_{A_i}$ and $I_{A_j}$ are comonotonic because $\{A_1, A_2, \cdots, A_n\}$ is a chain; to see that $f_p$ and $I_{A_i}$ are comonotonic, suppose not, then we can find $(x_1,x_2), (x_1',x_2')\in \mathbb{R}^2$ such that $f_p(x_1,x_2)\geq f_p(x_1',x_2')$ and $I_{A_i}(x_1,x_2)<I_{A_i}(x_1',x_2')$, the last inequality implies that $(x_1,x_2)\notin A_i$ and $(x_1',x_2')\in A_i$. This implies $f(x_1,x_2)<\alpha_i \leq f(x_1',x_2')$, and hence $f_p(x_1,x_2)<\alpha_i \leq f_p(x_1',x_2')$, a contradiction, so $f_p$ and $I_{A_i}$ are comonotonic.
Similarly, $f_p(\cdot,x_2)$ and $I_{A_i}(\cdot,x_2)$ are comonotonic for every $x_2\in \mathbb{R}$, for every $i=1, 2, \cdots, n$, i.e. $\{f_p(\cdot,x_2), I_{A_1}(\cdot,x_2), \cdots, I_{A_n}(\cdot,x_2)\}$ forms a comonotonic class, so does $\{f_p(x_1,\cdot), I_{A_1}(x_1,\cdot), \cdots, I_{A_n}(x_1,\cdot)\}$, $\forall x_1\in\mathbb{R}$.
%Because $\varphi_i$ are monotonic, so $f_p$ and $I_{A_i}$ are comonotonic, and then $f_p(\cdot,x_2)$ and $I_{A_i}(\cdot,x_2)$ are comonotonic for every $x_2\in\mathbb{R}$, for every $i=1,2,\cdots,n$. Hence $\{f_p, I_{A_1},\cdots,I_{A_n}\}$ forms a comonotonic class, so does $\{f_p(\cdot,x_2), I_{A_1}(\cdot,x_2),\cdots,I_{A_n}(\cdot,x_2)\}$.
By Lemma \ref{Lem1}, there are two probabilities $P$ and $Q$ on $\mathcal{F}$ such that
%$\{P_{x_2}\}_{x_2\in\mathbb{R}}$
%$$
%\mathbb{E}_{\mu}g(x,y)=\mathbb{E}_{P}g(x,y),
%$$
%similarly, there is a probability $Q$ on $\mathcal{A}_{Y}$ such that
$$
E_{V}[E_{V}[I_{ A_i}(x_1,X_2)]|_{x_1=X_1}] = E_{Q}[E_{P}[ I_{A_i}(x_1,X_2)]|_{x=X_1}],
$$
$i=1,2,\cdots,n$, and by comonotonic additivity, we have
$$
E_{V}[E_{V}[f_p(x_1,X_2)]|_{x_1=X_1}]=E_{Q}[E_{P}[f_p(x_1,X_2)]|_{x_1=X_1}].
$$
Define a capacity $\lambda$ on $\mathcal{U}_{f_p}$ as follows. For every $A\in\mathcal{U}_{f_p}=\{A_1,A_2,\cdots,A_n, \mathbb{R}^2, \emptyset\}$,
$$
\lambda(A)=E_{Q}[E_{P} [I_{A}(x_1,X_2)]|_{x_1=X_1}].
$$
%Hence, by Lemma \ref{Lem2} $\forall A\in\mathcal{U}_g$ is comonotonic, then $I_A$ is comonotonic, and then
%Because $A_{i+1}\subset A_{i}$, so $\{I_{A_{i}}\}_{i=1}^{n}$ are comonotonic, and then
By comonotonic additivity, we have
\begin{align*}
\int_{\mathbb{R}^2}f_p(x_1,x_2)d\lambda(x_1,x_2) &=\int_{\mathbb{R}^2} \left[\alpha_1 I_{A_1}(x_1,x_2) + \sum_{i=2}^{n}(\alpha_i -\alpha_{i-1})I_{A_i}(x_1,x_2)\right]d\lambda(x_1,x_2) \\
&=\alpha_1 \int_{\mathbb{R}^2} I_{A_1}(x_1,x_2)d\lambda(x_1,x_2) + \sum_{i=2}^{n}(\alpha_i -\alpha_{i-1})\int_{\mathbb{R}^2} I_{A_i}(x_1,x_2)d\lambda(x_1,x_2) \\
&=\alpha_1\lambda(A_1)+\sum_{i=2}^{n}(\alpha_i -\alpha_{i-1})\lambda(A_i) \\
&=\alpha_1 E_{Q}[E_{P}[ I_{A_1}(x_1,X_2)]|_{x_1=X_1}] + \sum_{i=2}^{n}(\alpha_i -\alpha_{i-1})E_{Q}[E_{P}[ I_{A_i}(x_1,X_2)]|_{x_1=X_1}] \\
&=E_{Q}\left[E_{P} \left[\alpha_1 I_{A_1}(x_1,X_2) +  \sum_{i=2}^{n}(\alpha_i -\alpha_{i-1})I_{A_i}(x_1,X_2) \right]\big|_{x_1=X_1}\right] \\
&=E_{Q}[E_{P}[f_p(x_1,X_2)]|_{x_1=X_1}]
=E_{V}[E_{V}[f_p(x_1,X_2)]|_{x_1=X_1}].
\end{align*}
%$$
%\mathbb{E}_{V}[\mathbb{E}_{V}[I_{A_i}(X_1,x)]|_{x=X_2}] = \mathbb{E}_{Q}[\mathbb{E}_{P}[ I_{A_i}(X_1,x)]|_{x=X_2}] = \lambda (A_i),
%$$
Because $X_1$ and $X_2$ are Fubini independent for $\mathcal{U}_{\hat{f}}$ and $A_i\in\mathcal{U}_{\hat{f}}$, so
\begin{equation}\label{Fubiniduli}
V((X_1,X_2)\in A)=\mathbb{E}_V[V((x_1,X_2)\in A)|_{x_1=X_1}]=\lambda(A),
\end{equation}
for all $A\in\mathcal{U}_{f_p}$. Hence
\begin{align*}
E_{V}[f_p(X_1,X_2)] &=E_V[\alpha_1I_{A_1}(X_1,X_2)+\sum_{i=2}^n(\alpha_i-\alpha_{i-1})I_{A_i}(X_1,X_2)] \\
&=\alpha_1 E_V[I_{A_1}(X_1,X_2)]+\sum_{i=2}^n(\alpha_i-\alpha_{i-1}) E_V[I_{A_i}(X_1,X_2)] \\
&=\alpha_1V((X_1,X_2)\in A_1) + \sum_{i=2}^n(\alpha_i-\alpha_{i-1})V((X_1,X_2)\in A_i) \\
&=\alpha_1\int_{\mathbb{R}^2}I_{A_1}(x_1,x_2)d\lambda(x_1,x_2) + \sum_{i=2}^n(\alpha_i-\alpha_{i-1})\int_{\mathbb{R}^2}I_{A_i}(x_1,x_2)d\lambda(x_1,x_2) \\
&=\int_{\mathbb{R}^2} f_p(x_1,x_2)d\lambda(x_1,x_2),
\end{align*}
by comonotonic additivity, hence we have %$\mathbb{E}_{V}[f_p(X_1,X_2)]=\mathbb{E}_{V}[\mathbb{E}_{V}[f_p(X_1,x)]|_{x=X_2}]$.
\begin{equation}\label{shizi3}
E_{V} [f_p(X_1,X_2)]=E_{V}[E_{V}[ f_p(x_1,X_2)]|_{x_1=X_1}].
\end{equation}

Following the Lemma 6.2 in Denneberg \cite{Den}, we have $\{f_p\}_{p=1}^{\infty}$ is an increasing sequence of measurable simple functions, which is the form as (\ref{shizi2}), converging uniformly to $\hat{f}$.
\iffalse
Because $\varphi_i, i=1,2$ are monotonic, so $f$ is comonotonic, and then each $f_p$ is comonotonic. In fact, one can easily see that, since $u_p$ is a monotonic transformation, the fact that $f(x,\cdot)$ and $f(x',\cdot)$ are comonotonic immediately implies that $f_p(x,\cdot)$ and $f_p(x',\cdot)$ are comonotonic.
\fi
% Applying the result we just proved, we have that for every $p$,
%\begin{equation}\label{shizi3}
%\mathbb{E}_{V} [f_p(X_1,X_2)]=\mathbb{E}_{V}[\mathbb{E}_{V}[ f_p(X_1,x)]|_{x=X_2}].
%\end{equation}
Now, for $x\in \mathbb{R}$ let
$$
G(x)=E_{V}[ \hat{f}(X_1,x)] \ \ and \ \ G_p(x)=E_{V}[ f_p(X_1,x)],
$$
and notice that $G_p(x)$ converges uniformly to $G(x)$. In fact, by the definition of $u_p$, for all $x_2\in \mathbb{R}$ we have
$$
\hat{f}(\cdot,x_2)-\frac{1}{2^p}\leq f_p(\cdot,x_2)\leq \hat{f}(\cdot,x_2)
$$
which implies
$$
G(x)-\frac{1}{2^p}\leq E_{V}[ f_p(X_1,x)] \leq G(x).
$$
This implies that $G$ is $\mathcal{B}$-measurable and
\begin{equation}\label{jixian1}
E_{V}[E_{V}[ f_p(x_1,X_2)]|_{x_1=X_1}] = E_{V}[ G_p(x_1)|_{x_1=X_1}]\rightarrow E_{V}[ G(x_1)|_{x_1=X_1}]=E_{V}[E_{V}[ \hat{f}(x_1,X_2)]|_{x_1=X_1}],
\end{equation}
as $p\rightarrow\infty$.

Because $f_p$ converges to $\hat{f}$ uniformly on $\mathbb{R}^2$, in fact for every $(x_1,x_2)\in \mathbb{R}^2$ we have
$$
\hat{f}(x_1,x_2)-\frac{1}{2^p}\leq f_p(x_1,x_2)\leq \hat{f}(x_1,x_2).
$$
This implies that
$$
E_{V} [\hat{f}(X_1,X_2)]-\frac{1}{2^p}\leq E_{V}[ f_p(X_1,X_2)]\leq E_{V}[ \hat{f}(X_1,X_2)],
$$
so that we get
\begin{equation}\label{jixian2}
E_{V}[ f_p(X_1,X_2)]\rightarrow E_{V}[\hat{f}(X_1,X_2)],
\end{equation}
as $p\rightarrow\infty$. Taking the limit as $p\rightarrow\infty$ on both sides of (\ref{shizi3}), combined with (\ref{jixian1}) and (\ref{jixian2}), we have completed the proof.
\\

\end{document}